\documentclass[a4paper,11pt,twoside,reqno]{amsart}

\usepackage[utf8]{inputenc}
\usepackage[plainpages=false,pdfpagelabels=true]{hyperref}
\usepackage{amssymb,amsthm}
\usepackage[margin=1in]{geometry}
\usepackage{slashed}

\newtheorem{Satz}{Theorem}[section]
\newtheorem{Prop}[Satz]{Proposition}
\newtheorem{Lem}[Satz]{Lemma}

\newtheorem{Cor}[Satz]{Corollary}
\theoremstyle{definition}

\newtheorem{Bem}[Satz]{Remark}
\newcommand{\tr}{\operatorname{Tr}}
\newcommand{\sff}{\mathrm{I\!I}}

\newcommand{\dv}{\text{ }dv_g}
\allowdisplaybreaks[1]

\renewcommand{\epsilon}{\varepsilon}

\newcommand{\N}{\ensuremath{\mathbb{N}}}
\newcommand{\s}{\mathbb{S}}

\numberwithin{equation}{section}

\newtheorem{innercustomgeneric}{\customgenericname}
\providecommand{\customgenericname}{}
\newcommand{\newcustomtheorem}[2]{%
  \newenvironment{#1}[1]
  {%
   \renewcommand\customgenericname{#2}%
   \renewcommand\theinnercustomgeneric{##1}%
   \innercustomgeneric
  }
  {\endinnercustomgeneric}
}

\newcustomtheorem{customthm}{Theorem}

\title{On the normal stability of the 4-harmonic and the ES-4-harmonic hypersphere}
\author{Volker Branding}
\date{\today}
\address{University of Vienna, Faculty of Mathematics\\
Oskar-Morgenstern-Platz 1, 1090 Vienna, Austria\\}
\email{volker.branding@univie.ac.at}
\subjclass[2010]{58E20; 53C43}
\keywords{4-harmonic maps; ES-4-harmonic maps; hypersphere; normal stability}
\thanks{The author gratefully acknowledges the support of the Austrian Science Fund (FWF) through the project "Geometric Analysis of Biwave Maps" (DOI: 10.55776/P34853)
}
\begin{document}

\begin{abstract}
Both 4-harmonic and ES-4-harmonic maps are two higher order generalizations of 
the well-studied
harmonic map equation given by a nonlinear elliptic partial differential equation
of order eight. Due to the large number of derivatives it is very difficult to 
find any difference in the qualitative behavior of these two variational problems.
It is well known that the small hypersphere \(\iota\colon\s^m(\frac{1}{2})\to\s^{m+1}\) 
is a critical point of both the 4-energy as well as the ES-4-energy but up to now
it has not been investigated if there is a difference concerning its stability.

The main contribution of this article is to show that the small hypersphere is unstable
with respect to normal variations both as 4-harmonic hypersphere
as well as ES-4-harmonic hypersphere and that its normal index equals one in both cases.
\end{abstract} 

\maketitle

\section{Introduction and results}
At the heart of the geometric calculus of variations is the aim to find
interesting maps between Riemannian manifolds, one way of achieving this
is to extremize a given energy functional.

For maps between Riemannian manifolds \(\phi\colon (M,g)\to (N,h)\) a lot of attention has been paid
to the energy of a map 
\begin{align}
\label{energy-map}
E(\phi)=\frac{1}{2}\int_M|d\phi|^2\dv.
\end{align}
The critical points of \eqref{energy-map} are characterized by the vanishing
of the so-called \emph{tension field} which is defined by
\begin{align}
\label{tension-field}
0=\tau(\phi):=\tr_g\bar\nabla d\phi,
\end{align}
where \(\bar\nabla\) represents the connection on the pull-back bundle \(\phi^\ast TN\).
Solutions of \eqref{tension-field} are called \emph{harmonic maps}
and the latter have been studied intensively in the literature.
The harmonic map equation is a second order semilinear elliptic partial differential
equation. Due to its second order nature powerful tool such as the maximum principle 
are available for its analysis. The famous  result of Eells - Sampson \cite{MR164306} guarantees the existence of a harmonic map in each homotopy class of maps if both \(M,N\) are closed and if \(N\)
has non-positive sectional curvature. Note that these harmonic maps are stable in the sense
that the second variation of the energy \eqref{energy-map} will always be positive.
For more details and references on the stability of harmonic maps one can consult \cite{MR4477489}.
For an overview on the current status of research on harmonic maps
we refer to  \cite{MR2044031} and \cite{MR2389639}.

Recently, many researchers got attracted in energy functionals that contain higher
derivatives generalizing the energy of a map \eqref{energy-map}.
In this regard, much attention has been paid to the \emph{bienergy} of a map
\begin{align}
E_{2}(\phi)=\frac{1}{2}\int_M|\tau(\phi)|^2\dv,    
\end{align}
whose critical points are characterized by a fourth order equation and which are called
\emph{biharmonic maps}.
A direct application of the maximum principle shows that biharmonic maps from compact Riemannian manifolds to targets of non-positive curvature have to be harmonic \cite{MR886529}. However, in the case of a target with positive curvature, i.e. a Euclidean sphere, there are additional non-harmonic solutions, 
see for example \cite{MR4410183}.
The stability of biharmonic maps to spheres has been studied intensively, see for example the article \cite{MR2135286}.
For more details and the current status of research on biharmonic maps in Riemannian geometry we refer to the book \cite{MR4265170} and the survey \cite{MR4410183}.

A general higher order version of harmonic maps is given by the so-called
\emph{polyharmonic maps of order k} or just \emph{k-harmonic maps}. 
These are critical points of the following energy functionals,
where we need to distinguish between polyharmonic maps of even and odd order.
In the even case \((k=2s,s\in\N)\) we set
\begin{align}
\label{poly-energy-even}
E_{2s}(\phi)=\frac{1}{2}\int_M|\bar\Delta^{s-1}\tau(\phi)|^2\dv,
\end{align}
whereas in the odd case \((k=2s+1,s\in\N)\) we have
\begin{align}
\label{poly-energy-odd}
E_{2s+1}(\phi)=\frac{1}{2}\int_M|\bar\nabla\bar\Delta^{s-1}\tau(\phi)|^2\dv.
\end{align}
Here, we use \(\bar\Delta\) to denote the connection Laplacian on the pull-back bundle \(\phi^\ast TN\).
The first variation of \eqref{poly-energy-even}, \eqref{poly-energy-odd} was calculated in \cite{MR3007953}. For the precise form of the critical points of
\eqref{poly-energy-even}, \eqref{poly-energy-odd}
and the current status of research on polyharmonic maps
we refer to \cite{MR4106647}.

For explicit solutions to the \(k\)-harmonic map equation we refer to
\cite{annarharm,MR4552081,MR4730425,MR3711937,MR3790367,MR4462636}.

Another possible generalization of harmonic maps,
first suggested by Eells and Sampson in 1964 \cite{MR172310},
can be obtained by studying the critical points of the following energy functional
\begin{align}
\label{k-energy}
E^{ES}_k(\phi):=\frac{1}{2}\int_M|(d+d^\ast)^k\phi|^2\dv
=\frac{1}{2}\int_M|(d+d^\ast)^{k-2}\tau(\phi)|^2\dv,\qquad k=1,2,\ldots.
\end{align}
The critical points of \eqref{k-energy} are referred to as \emph{ES-k-harmonic maps}.

In the following we will provide a brief overview on the structure of \eqref{k-energy}.
For \(k=1\) the energy \eqref{k-energy} reduces to the energy of a map \eqref{energy-map}.

In the case of \(k=2\), which is also obtained in \eqref{poly-energy-even} for \(s=1\),
we are led to the aforementioned bienergy \(E_2=E^{ES}_2\).

For \(k=3\) we gain the trienergy of a map \(E_3=E^{ES}_3\), which corresponds to \eqref{poly-energy-odd}
with \(s=1\), and its critical points are called \emph{triharmonic maps}.
An overview on triharmonic maps is provided by \cite[Section 4]{MR4007262}, \cite{MR4634317} and references therein.

However, for \(k\geq 4\) the energy functional \eqref{k-energy} contains additional curvature terms
and can in general no longer be written in the form \eqref{poly-energy-even}, \eqref{poly-energy-odd}, that is \(E_r\neq E_r^{ES}\) for \(r\geq 4\).
An extensive analysis of \eqref{k-energy} and its critical points was carried out in \cite{MR4106647}.

Note that in the case of \(\dim M=1\) the functionals \eqref{poly-energy-even}, \eqref{poly-energy-odd} and \eqref{k-energy} coincide and the critical points are given by polyharmonic curves.
For the current status of research on this subject we refer to the recent article
\cite{MR4542687} and references cited in that manuscript.

It can directly be seen that every harmonic map, that is a solution of \(\tau(\phi)=0\),
automatically is a critical point of \eqref{poly-energy-even}, \eqref{poly-energy-odd}
as well as of \eqref{k-energy}. Hence, we are interested in determining those critical 
points which are non-harmonic, the latter are usually referred to as \emph{proper k-harmonic}
or \emph{proper ES-k-harmonic}, respectively.

This article is devoted to polyharmonic maps of order \(4\) arising either as critical points of \eqref{poly-energy-even} for \(s=2\)
or as critical points of \eqref{k-energy} for \(k=4\).

The energy functional for 4-harmonic maps (corresponding to \eqref{poly-energy-even} with \(s=2\)) is given by
\begin{align}
\label{energy-4-harmonic}
E_4(\phi)=\frac{1}{2}\int_M|\bar\Delta\tau(\phi)|^2\dv.
\end{align}

The critical points of \eqref{energy-4-harmonic} are characterized by the 
vanishing of the 4-tension field
\begin{align}
\label{4-tension}
0=\tau_{4}(\phi):=&\bar\Delta^{3}\tau(\phi)-R^N(\bar\Delta^{2}\tau(\phi),d\phi(e_j))d\phi(e_j) \\
\nonumber&+R^N(\tau(\phi),\bar\nabla_{e_j}\bar\Delta\tau(\phi))d\phi(e_j) 
-R^N(\bar\nabla_{e_j}\tau(\phi),\bar\Delta\tau(\phi))d\phi(e_j),
\end{align}
where \(R^N\) represents the Riemann curvature tensor of the manifold \(N\)
and \(\{e_i\},i=1,\ldots,\dim M\) is an orthonormal frame field tangent to \(M\).
Solutions of \eqref{4-tension} are called \emph{4-harmonic maps}.
For a derivation of the 4-tension field \eqref{4-tension} we refer to \cite[Theorem 2.5]{MR3007953}.

The energy functional for ES-4-harmonic maps (corresponding to \eqref{k-energy} with \(k=4\)) is given by
\begin{align}
\label{energy-es4-harmonic}
E^{ES}_4(\phi)&=\frac{1}{2}\int_M|(d+d^\ast)^4\phi|^2\dv\\
\nonumber&=\underbrace{\frac{1}{2}\int_M|\bar\Delta\tau(\phi)|^2\dv}_{:=E_4(\phi)}
+\underbrace{\frac{1}{4}\int_M|R^N(d\phi(e_i),d\phi(e_j))\tau(\phi)|^2\dv}_{:=\hat E_4(\phi)}.
\end{align}

The first variation of \eqref{energy-es4-harmonic} 
was calculated in \cite[Section 3]{MR4106647} and is characterized by the vanishing 
of the ES-4-tension field \(\tau_4^{ES}(\phi)\) given by the following expression
\begin{align}
\label{es-4-tension}
\tau_4^{ES}(\phi)=
\tau_4(\phi)+\hat\tau_4(\phi).
\end{align}
Here, \(\tau_4(\phi)\) denotes the 4-tension field \eqref{4-tension} and the quantity \(\hat\tau_4(\phi)\)
is defined by
\begin{align*}
\hat{\tau}_4(\phi)=-\frac{1}{2}\big(2\xi_1+2d^\ast\Omega_1+\bar\Delta\Omega_0+\tr R^N(d\phi(\cdot),\Omega_0)d\phi(\cdot)\big),
\end{align*}
where we have used the following abbreviations
\begin{align*}
\Omega_0&=R^N(d\phi(e_i),d\phi(e_j))R^N(d\phi(e_i),d\phi(e_j))\tau(\phi), \\
\nonumber\Omega_1(X)&=R^N(R^N(d\phi(X),d\phi(e_j))\tau(\phi),\tau(\phi))d\phi(e_j),\\
\nonumber\xi_1&=-(\nabla_{d\phi(e_j)}R^N)(R^N(d\phi(e_i),d\phi(e_j))\tau(\phi),\tau(\phi))d\phi(e_i).
\end{align*}

Up to now, not many properties of \(ES-4\)-harmonic maps are known in the mathematics literature.
The stress-energy tensor associated with \(ES-4\)-harmonic maps was derived in \cite{MR4293944} 
and with its help the structure of finite energy solutions from Euclidean space has been investigated.
In particular, it was shown that one needs to impose different energy conditions
that force a \(4\)-harmonic map to be harmonic compared to the case of \(ES-4\)-harmonic maps.

In \cite{MR4550874} it was shown that both \(4\)-harmonic as well as \(ES-4\)-harmonic
maps satisfy the unique continuation property.

The general structure of r-harmonic and ES-r-harmonic hypersurfaces in space forms
has been investigated in \cite{MR4462636}.

Concerning the existence of \(4\)-harmonic and \(ES-4\)-harmonic maps we have
\begin{Satz}{\cite[Theorem 2.1]{MR4106647}, \cite[Theorem 1.1]{MR3711937}}
The small hypersphere \(\iota\colon\s^m(\frac{1}{2})\to\s^{m+1}\)
is both a proper \(4\)-harmonic as well as a proper \(ES-4\)-harmonic submanifold 
of \(\s^{m+1}\).
\end{Satz}

Note that this hypersphere is totally umbilic and thus has constant mean curvature
which we will abbreviate by CMC throughout this article.

Although the above mentioned small hypersphere is a critical point of both
\(E_4\) and \(E^{ES}_4\) one could expect that there is a difference concerning
its stability. In this manuscript we will systematically approach this question by deriving
the second variation of both energy functionals and in a second step we will then
study the stability with respect to normal variations of the small hypersphere
considered as a critical point of either \(E_4\) or \(E^{ES}_4\).

The stability of a given 4-harmonic map is characterized by the second variation of the 4-energy of a map \eqref{energy-4-harmonic} evaluated at a critical point
which we denote by 
\(Q_4(V,W)\), where
\(V,W\in\Gamma(\phi^\ast TN)\).
We say that a 4-harmonic map \(\phi\) is \emph{stable} if
\begin{align*}
Q_4(V,V)>0\qquad \textrm{ for all } V\in\Gamma(\phi^\ast TN),\qquad V\neq 0
\end{align*}
and \emph{weakly stable} if
\begin{align*}
Q_4(V,V)\geq 0\qquad \textrm{ for all } V\in\Gamma(\phi^\ast TN),\qquad V\neq 0.
\end{align*}
We will use the same terminology if we consider the energy \(E^{ES}_4\)
instead of \(E_4\).

We will now present the main results obtained in this article.
Our first result is given by the following
\begin{customthm}{\ref{thm:harmonic-stable-es4}}
A harmonic map \(\phi\colon M\to N\) is a weakly stable critical point of both 
\(E_4(\phi)\) and \(E^{ES}_4(\phi)\).
\end{customthm}

Furthermore, we also show that a harmonic map is a weakly stable critical point
of \(\hat E_4(\phi)\), see Theorem \ref{thm:harmonic-stable-escurv} for the precise details.

In order to present our main results on the stability of the small hypersphere 
we recall the concept of \emph{normal stability}.
Following the terminology used for minimal, biharmonic \cite{MR4386842}
and triharmonic \cite{MR4634317} hypersurfaces
we define the normal index of a proper 4-harmonic hypersurface
as well as an ES-4-harmonic hypersurface
to be the maximal dimension of any linear subspace 
on which the second variation is negative, that is
\begin{align}
\operatorname{Ind}_4^{\rm{nor}}(M\to N)
:=\max\{\dim L, L\subset C^\infty_0(M)\mid Q_4(f\nu,f\nu)<0
,~~\forall f\in L\},
\end{align}
where \(\nu\) represents the normal of the hypersurface. In addition, we use 
\(\operatorname{Ind}_{ES-4}^{\rm{nor}}(M\to N)\) to represent the normal index
of \(ES-4\)-harmonic hypersurfaces, respectively.

Let us now give the two main results of this article
which show that the normal index of the small hypersphere \(\iota\colon\s^m(\frac{1}{2})\to\s^{m+1}\) is equal to one both when considered as a critical point
of the 4-energy \eqref{energy-4-harmonic} or as a critical point
of the \(ES-4\)-energy \eqref{energy-es4-harmonic}.
More precisely, we have
\begin{customthm}{\ref{thm:4hyper-stab}}
Let \(\phi\colon \s^m(1/2)\hookrightarrow\s^{m+1}\)
be the small proper 4-harmonic hypersphere.
The index characterizing its normal stability is equal to one, i.e.
\begin{align*}
\operatorname{Ind}_4^{nor}\big(\s^m(1/2)\hookrightarrow\s^{m+1}\big)=1.
\end{align*}
\end{customthm}

\begin{customthm}{\ref{thm:stability-es4}}
Let \(\phi\colon \s^m(1/2)\hookrightarrow\s^{m+1}\)
be the small proper ES-4-harmonic hypersphere.
The index characterizing its normal stability is equal to one, i.e.
\begin{align*}
\operatorname{Ind}_{ES-4}^{nor}\big(\s^m(1/2)\hookrightarrow\s^{m+1}\big)=1.
\end{align*}
\end{customthm}

In addition, we also show that
the small hypersphere \(\phi\colon \s^m(1/2)\hookrightarrow\s^{m+1}\)
is a weakly stable critical point of \(\hat E_4(\phi)\) with respect to normal variations,
see Theorem \ref{thm:hypersphere-crit-curv} for the precise details.

\begin{Bem}
The normal stability of biharmonic hyperspheres was recently studied in 
\cite{MR4386842}, see also \cite{MR2135286} for a general analysis of the
stability of biharmonic maps in spheres.

Concerning the normal stability of triharmonic hypersurfaces in space forms 
we refer to the recent article \cite{MR4634317}.
In particular, it was conjectured in \cite[Remark 1.1]{MR4634317} that
the normal index of the small \(k\)-harmonic hypersphere 
\(\phi\colon\s^m({1/\sqrt{k}})\to\s^{m+1}\) is always equal to one 
and the main results of this article further support this conjecture.

\end{Bem}

Throughout this article we will use the following notation:
Indices on the domain manifold will be denoted by Latin letters \(i=1,\ldots,m=\dim M\).
We will use the following sign convention for the \emph{rough Laplacian} acting on sections of $\phi^{\ast}TN$
\begin{align*}
\bar{\Delta}=d^\ast d =-\big(\bar\nabla_{e_i}\bar\nabla_{e_i}-\bar\nabla_{\nabla^M_{e_i}e_i}\big),
\end{align*}
where $\{e_i\},i=1,\ldots m$ is a local orthonormal frame field tangent to $M$.
Moreover, we employ the summation convention and tacitly sum over repeated indices.
Note that, due to our choice of sign convention, the Laplace operator acting
on functions has a positive spectrum.
Throughout this article we make use of the following sign convention for the Riemann curvature tensor
\begin{align*}
R(X,Y)Z=\nabla_X\nabla_YZ-\nabla_Y\nabla_XZ-\nabla_{[X,Y]}Z
\end{align*}
for given vector fields \(X,Y,Z\).

This article is organized as follows: In Section 2 we collect a number of preliminary results (commutation formulas, geometry of hypersurfaces). 
Section 3 is devoted to the second variation of the 4-energy while Section 4 
derives the corresponding results for the ES-4-energy.

\section{Preliminaries}
In this section we will collect and establish various results which will be
frequently employed in the rest of the manuscript.

For most of our statements we will assume that the target manifold \(N\) is a space
form of constant curvature \(K\).
In this case the Riemann curvature tensor acquires the simple form
\begin{align}
\label{curvature-space-form}
R^N(X,Y)Z=K(\langle Y,Z\rangle X-\langle X,Z\rangle Y),
\end{align}
where \(X,Y,Z\) are vector fields on \(N\) and \(K\)
represents the constant curvature of the space form \(N\).

We will often consider a variation of a map \(\phi\colon M\to N\) which is
a map \(\phi_t\colon (-\epsilon,\epsilon)\times M\to N,\epsilon>0\) together with
its variational vector field
\begin{align}
\label{dfn:variation-phi}
\frac{\nabla\phi_t}{\partial t}\big|_{t=0}=V,
\end{align}
where \(V\in\Gamma(\phi^\ast TN)\).
Throughout this article \(\{e_i\},i=1,\ldots, M\) 
will always represent a local orthonormal frame field tangent to \(M\)
at a fixed but arbitrary point \(p\in M\) such that at this point we have
\begin{align*}
\nabla_{e_i}e_j=0,\qquad i,j=1,\ldots,\dim M.
\end{align*}

\subsection{Commutator formulas}
In the following we will often make use of the following well-known
formula
\begin{align}
\label{commutator-t-tension}
\frac{\bar\nabla}{\partial t}\tau(\phi_t)=&-\bar\Delta d\phi_t(\partial_t)
+R^N(d\phi_t(\partial_t),d\phi_t(e_i))d\phi_t(e_i).
\end{align}

A direct consequence of this commutator formula is the following
\begin{Lem}
Consider a variation of the map \(\phi\colon M\to N\) as defined in
\eqref{dfn:variation-phi}. Then, the following formulas hold
\begin{align}
\label{commutator-second-variation}
\frac{\bar\nabla}{\partial t}\bar\nabla_{X}\tau(\phi_t)=&
-\bar\nabla_{X}\bar\Delta d\phi_t(\partial_t)
+\bar\nabla_{X}\big(R^N(d\phi_t(\partial_t),d\phi_t(e_k))d\phi_t(e_k)\big)
+R^N(d\phi_t(\partial_t),d\phi_t(X))\tau(\phi_t), \\
\nonumber\frac{\bar\nabla}{\partial t}\bar\Delta\tau(\phi_t)=&
-R^N(d\phi_t(\partial_t),d\phi_t(e_j))\bar\nabla_{e_j}\tau(\phi_t)
-\bar\nabla_{e_j}\big(R^N(d\phi_t(\partial_t),d\phi_t(e_j))\tau(\phi_t)\big) \\
\nonumber&+\bar\Delta\big(R^N(d\phi_t(\partial_t),d\phi_t(e_j))d\phi_t(e_j)\big)
-\bar\Delta^2d\phi_t(\partial_t), \\
\nonumber\frac{\bar\nabla}{\partial t}\bar\Delta^k\tau(\phi_t)=&
-R^N(d\phi_t(\partial_t),d\phi_t(e_j))\bar\nabla_{e_j}\bar\Delta^{k-1}\tau(\phi_t)
-\bar\nabla_{e_j}\big(R^N(d\phi_t(\partial_t),d\phi_t(e_j))\bar\Delta^{k-1}\tau(\phi_t)\big) \\
&+\nonumber\bar\Delta\frac{\tilde\nabla}{\partial t}\bar\Delta^{k-1}\tau(\phi_t)
\end{align}
for all \(X\in TM\) and \(k\in\N\) with \(k\geq 1\).
\end{Lem}
\begin{proof}
This follows from \eqref{commutator-t-tension} and arguments similar to
the ones used in the proof of \cite[Lemma 2.2]{MR4634317}.
\end{proof}

\subsection{Some facts on hypersurfaces}
Let us recall a number of geometric facts on hypersurfaces \(M^m\) in a Riemannian manifold 
\(N^{m+1}\) mostly in order to fix the notation.
The connections on \(M^m\) and \(N^{m+1}\) are related by the following formula
\begin{align}
\label{dfn:sff}
\nabla^N_XY=\nabla^M_XY+\sff(X,Y),
\end{align}
where \(\sff\) represents the second fundamental form of the hypersurface
and \(X,Y\) are vector fields on \(M^m\).

Let \(\nu\) be the global unit normal of the hypersurface \(M^m\), then 
its shape operator \(A\) is given by
\begin{align}
\label{dfn:shape}
\nabla_{X}\nu=-A(X),
\end{align}
where \(X\) is a vector field on \(M\).

The shape operator \(A\) and the second fundamental form \(\sff\)
are related by
\begin{align}
\label{relation-shape-sff}
\sff(X,Y)=\langle A(X),Y\rangle\nu.
\end{align}

If \(\phi\colon M^{m}\to N^{m+1}\) is an isometric immersion
the tension field acquires the form \(\tau(\phi)=mH\nu\),
where \(H\) represents the mean curvature function of the hypersurface.
In addition, we can compute the mean curvature function by \(H=\frac{1}{m}\tr A\).

In the following lemma we provide without proof three standard facts which we shall use in this section, 
see for example \cite[Lemma 2.1]{MR4462636}
or \cite[Lemma 4.1]{MR4552081}.

\begin{Lem}
\label{lemma-sff}
Let $\phi\colon M^m\to N^{m+1}$ be a hypersurface in a Riemannian manifold \(N^{m+1}\).
In addition, let $A$ be the shape operator and $H=1/m\tr A$ the mean curvature function. Then, we have that 
\begin{enumerate}
\item $(\nabla A) (\cdot,\cdot)$ is symmetric,
\item $\langle (\nabla A) (\cdot,\cdot), \cdot \rangle$ is totally symmetric,
\item $\tr (\nabla A) (\cdot,\cdot)= m\operatorname{grad} H$. 
\end{enumerate}
\end{Lem}

This allows us to prove the following

\begin{Lem}
Let $\phi\colon M^m\to N^{m+1}$ be a hypersurface. 
Suppose that \(k\in C^\infty(M)\) and let \(\nu\) be the unit normal
of the hypersurface. Then, the following identities hold
\begin{align}
\label{eq:hyper-general-identities}
\bar\Delta(k\nu)=&(\Delta k+k|A|^2)\nu+2A(\operatorname{grad} k)+mk\operatorname{grad} H,\\
\nonumber\bar\nabla_{X}\bar\Delta(k\nu)=&
(\nabla_{X}\Delta k+(\nabla_{X} k)|A|^2+k\nabla_{X}|A|^2)\nu 
-(\Delta k+k|A|^2)A(X) \\
\nonumber&+2(\nabla_{X}A)(\operatorname{grad} k)
+2A(\nabla_{X}\operatorname{grad} k)
+2\sff(X,A(\operatorname{grad} k))\\
&\nonumber+m(\nabla_{X}k)\operatorname{grad} H
+mk\nabla_{X}\operatorname{grad} H
+mk\sff(X,\operatorname{grad} H),\\
\nonumber\bar\Delta^2(k\nu)=&\big(\Delta^2 k+|A|^2\Delta k+\Delta (|A|^2k)+|A|^4k\big)\nu
\\
\nonumber &+2A(\operatorname{grad}\Delta k)+m\Delta k\operatorname{grad} H
+2A(\operatorname{grad}(|A|^2k))+m|A|^2k\operatorname{grad}H
\\
\nonumber&
-2(\nabla_{e_i}\nabla_{e_i}A)(\operatorname{grad}k)
-4(\nabla_{e_i}A)(\nabla_{e_i}\operatorname{grad}k)
-4\langle A(e_i),(\nabla_{e_i}A)(\operatorname{grad}k)\rangle\nu \\
\nonumber&+2A(\Delta\operatorname{grad}k)
-4\langle A(e_i),A(\nabla_{e_i}\operatorname{grad}k)\rangle\nu
-2m\operatorname{grad}HA(\operatorname{grad}k)\nu\\
\nonumber&
+2\langle A(e_i),A(\operatorname{grad}k)\rangle A(e_i) 
+m\bar\Delta\big(k\operatorname{grad} H\big),
\end{align}
where \(X\) represents a vector field on \(M\).
\end{Lem}
\begin{proof}
The first two statements were already obtained in \cite[Lemma 3.2]{MR4634317}.
To obtain the third statement we apply the first formula of \eqref{eq:hyper-general-identities} and find
\begin{align*}
\bar\Delta^2(k\nu)=\bar\Delta(\Delta k\nu)+\bar\Delta(k|A|^2\nu)+2\bar\Delta\big(A(\operatorname{grad} k)\big)
+m\bar\Delta\big(k\operatorname{grad} H\big).
\end{align*}
Again, applying the first identity of \eqref{eq:hyper-general-identities} once more we find
\begin{align*}
\bar\Delta(\Delta k\nu)=&(\Delta^2 k+|A|^2\Delta k)\nu+2A(\operatorname{grad}\Delta k)+m\Delta k\operatorname{grad} H,
\\
\bar\Delta(k|A|^2\nu)=&\big(\Delta (|A|^2k)+|A|^4k\big)\nu+2A(\operatorname{grad}(|A|^2k))+m|A|^2k\operatorname{grad} H.
\end{align*}
In order to obtain the precise expression for \(\bar\Delta\big(A(\operatorname{grad} k)\big)\)
we first of all calculate
\begin{align*}
\bar\nabla_{e_i}\big(A(\operatorname{grad} k)\big)=
(\nabla_{e_i}A)(\operatorname{grad} k)
+\langle A(e_i),A(\operatorname{grad} k)\rangle\nu
+A(\nabla_{e_i}\operatorname{grad}k).
\end{align*}

Furthermore, a direct calculation yields
\begin{align*}
\bar\Delta\big(A(\operatorname{grad} k)\big)=&-(\nabla_{e_i}\nabla_{e_i}A)(\operatorname{grad}k)
-2(\nabla_{e_i}A)(\nabla_{e_i}\operatorname{grad}k)
-2\langle A(e_i),(\nabla_{e_i}A)(\operatorname{grad}k)\rangle\nu \\
&+A(\Delta\operatorname{grad}k)-2\langle A(e_i),A(\nabla_{e_i}\operatorname{grad}k)\rangle\nu
-m\operatorname{grad}HA(\operatorname{grad}k)\nu\\
&
+\langle A(e_i),A(\operatorname{grad}k)\rangle A(e_i)
\end{align*}
completing the proof.
\end{proof}

In the upcoming analysis of the normal stability of the \(4-\) and the \(ES-4\) harmonic
hyperspheres it suffices to consider the following simplified version 
of the previous Lemma.

\begin{Cor}
Let $\phi\colon M^m\to N^{m+1}$ be a hypersurface with parallel shape
operator, that is \(\nabla A=0\).
Suppose that \(k\in C^\infty(M)\) and let \(\nu\) be the unit normal
of the hypersurface. Then, the following identities hold
\begin{align}
\label{eq:laplace-hypersurface}
\bar\Delta(k\nu)=&(\Delta k+k|A|^2)\nu+2A(\operatorname{grad} k),\\
\label{eq:nabla-laplace-hypersurface}
\bar\nabla_{X}\bar\Delta(k\nu)=&
(\nabla_{X}\Delta k+(\nabla_{X} k)|A|^2)\nu 
-(\Delta k+k|A|^2)A(X) \\
\nonumber&+2A(\nabla_{X}\operatorname{grad} k)
+2\sff(X,A(\operatorname{grad} k)), \\
\label{eq:laplace-2-hypersurface}
\bar\Delta^2(k\nu)=&\big(\Delta^2 k+2|A|^2\Delta k+|A|^4k\big)\nu
-4\langle A(e_i),A(\nabla_{e_i}\operatorname{grad}k)\rangle\nu
\\
\nonumber &+2A(\operatorname{grad}\Delta k)
+2|A|^2A(\operatorname{grad}k)
+2A(\Delta\operatorname{grad}k) \\
\nonumber&
+2\langle A(e_i),A(\operatorname{grad}k)\rangle A(e_i).
\end{align}
\end{Cor}
\begin{proof}
By assumption the shape operator of the hypersurface is parallel, i.e. \(\nabla A=0\)
which directly implies that \(\operatorname{grad} H=0\) and \(|A|^2=const\).
Using these identities in the formulas provided by the previous Lemma already completes the proof.
\end{proof}

Moreover, let us also recall the following fact:
Consider the inclusion \(\iota\colon\s^m(a)\hookrightarrow\s^{m+1}\),
where \(0<a<1\),
then the corresponding shape operator is given by
\begin{align*}
A=-\frac{\sqrt{1-a^2}}{a}\operatorname{Id}.
\end{align*}

\subsection{Spectral properties of the Laplace-Beltrami operator}

After having established the necessary geometric formulas for our analysis
we also recall a number of results on the spectrum of the Laplace operator on closed manifolds.
On an arbitrary closed Riemannian manifold \((M,g)\)
its eigenvalues are discrete, have finite multiplicity and satisfy
\begin{align*}
0=\lambda_0\leq\lambda_1\leq\lambda_2\leq\ldots\to\infty.
\end{align*}

In the case of the \(p\)-dimensional sphere of radius \(R\), that is \(\s^p(R)\),
we know that the eigenvalues of the Laplacian are given by
\begin{align}
\label{ev-laplacian-sphere}
\lambda_j=\frac{j(j+p-1)}{R^2},\qquad j=0,1,2,\ldots
\end{align}
with multiplicities \(m_{\lambda_0}=1,m_{\lambda_1}=p+1\)
and
\begin{align*}
m_{\lambda_j}={{p+j}\choose{j}}-{{p+j-2}\choose{j-2}},\qquad j=2,\ldots.
\end{align*}

For the general properties of the spectrum of the Laplacian we refer to the book \cite{MR0282313}, for more details concerning the Laplace operator in the context of
the stability of CMC hypersurfaces one may consult
the survey article \cite{MR2301375}.

\section{The second variation of the 4-energy}
Let us start by recalling the first variation formula 
\eqref{4-tension}
for \(4\)-harmonic maps
\begin{align}
\label{first-variation-4-energy}
\frac{d}{dt}E_4(\phi_t)=-\int_M&\big(\langle d\phi_t(\partial_t),
\bar\Delta^{3}\tau(\phi_t)-R^N(\bar\Delta^{2}\tau(\phi_t),d\phi_t(e_j))d\phi_t(e_j) \\
\nonumber&+R^N(\tau(\phi_t),\bar\nabla_{e_j}\bar\Delta\tau(\phi_t))d\phi_t(e_j) 
-R^N(\bar\nabla_{e_j}\tau(\phi_t),\bar\Delta\tau(\phi_t))d\phi_t(e_j)\rangle\big)\dv,
\end{align}
where \(\phi_t\) is a variation of the map \(\phi\) as defined in \eqref{dfn:variation-phi}.

\subsection{Second variational formula - the general case}
In this subsection we derive the second variation formula of the 4-energy \eqref{energy-4-harmonic}
for a general target manifold.

\begin{Lem}[Second variation of the 4-energy (rough version)]
Let \(\phi\colon M\to N\) be a smooth \(4\)-harmonic map and consider a variation of the map \(\phi\)
as defined in \eqref{dfn:variation-phi}. Then, the second variation of \eqref{energy-4-harmonic}
is given by
\begin{align}
\label{eq:snd-variation-4-rough}
\frac{d^2}{dt^2}\big|_{t=0}E_4(\phi_t)=
&\int_M|\bar\Delta^2V|^2\dv\\
\nonumber&+\int_M\langle V,R^N(V,d\phi(e_j))\bar\nabla_{e_j}\bar\Delta^{2}\tau(\phi)
\rangle\dv\\
\nonumber&-\int_M\langle \bar\nabla_{e_j} V,R^N(V,d\phi(e_j))\bar\Delta^{2}\tau(\phi) 
\rangle\dv\\
\nonumber&+\int_M\langle \bar\Delta V,R^N(V,d\phi(e_j))\bar\nabla_{e_j}\bar\Delta\tau(\phi)
\rangle\dv\\
\nonumber&-\int_M\langle \bar\nabla_{e_j}\bar\Delta V,R^N(V,d\phi(e_j))\bar\Delta\tau(\phi)
\rangle\dv\\
\nonumber&+\int_M\langle \bar\Delta^2 V,R^N(V,d\phi(e_j))\bar\nabla_{e_j}\tau(\phi)
\rangle\dv\\
\nonumber&-\int_M\langle \bar\nabla_{e_j}\bar\Delta^2 V,R^N(V,d\phi(e_j))\tau(\phi)
\rangle\dv\\
\nonumber&-\int_M\langle \bar\Delta^3 V,R^N(V,d\phi(e_k))d\phi(e_k)
\rangle\dv\\
\nonumber&+\int_M\langle V,(\nabla_{V}R^N)(\bar\Delta^{2}\tau(\phi),d\phi(e_j))d\phi(e_j)
\rangle\dv \\
\nonumber&-\int_M\langle V,R^N(R^N(V,d\phi(e_k))\bar\nabla_{e_k}\bar\Delta\tau(\phi),d\phi(e_j))d\phi(e_j)
\rangle\dv \\
\nonumber&-\int_M\langle V,R^N(\bar\nabla_{e_k}\big(R^N(V,d\phi(e_k))\bar\Delta\tau(\phi)\big),d\phi(e_j))d\phi(e_j)
\rangle\dv \\
\nonumber&-\int_M\langle V,R^N(\bar\Delta\big(R^N(V,d\phi(e_k))\bar\nabla_{e_k}\tau(\phi)\big),d\phi(e_j))d\phi(e_j)
\rangle\dv \\
\nonumber&-\int_M\langle V,R^N(\bar\Delta\bar\nabla_{e_k}\big(R^N(V,d\phi(e_k))\tau(\phi)\big),d\phi(e_j))d\phi(e_j)
\rangle\dv \\
\nonumber&-\int_M\langle V,R^N(\bar\Delta^3V,d\phi(e_j))d\phi(e_j)
\rangle\dv \\
\nonumber&+\int_M\langle V,R^N(\bar\Delta^{2}\big(R^N(V,d\phi(e_k))d\phi(e_k)\big),d\phi(e_j))d\phi(e_j)
\rangle\dv \\
\nonumber&
+\int_M\langle V,R^N(\bar\Delta^2\tau(\phi),\bar\nabla_{e_j}V)d\phi(e_j)
\rangle\dv \\
\nonumber&
+\int_M\langle V,R^N(\bar\Delta^2\tau(\phi),d\phi(e_j))\bar\nabla_{e_j}V
\rangle\dv\\
\nonumber&
-\int_M\langle V,(\nabla_{V}R^N)(\tau(\phi),\bar\nabla_{e_j}\bar\Delta\tau(\phi))d\phi(e_j)
\rangle\dv\\
\nonumber&+\int_M\langle V,R^N(\bar\Delta V,\bar\nabla_{e_j}\bar\Delta\tau(\phi))d\phi(e_j)
\rangle\dv\\
\nonumber&-\int_M\langle V,R^N(R^N(V,d\phi(e_k))d\phi(e_k),\bar\nabla_{e_j}\bar\Delta\tau(\phi))d\phi(e_j)
\rangle\dv\\
\nonumber&-\int_M\langle V,R^N(\tau(\phi),R^N(V,d\phi(e_j))\bar\Delta\tau(\phi))d\phi(e_j)
\rangle\dv\\
\nonumber&+\int_M\langle V,R^N(\tau(\phi),\bar\nabla_{e_j}\big(R^N(V,d\phi(e_k))\bar\nabla_{e_k}\tau(\phi)\big))d\phi(e_j)
\rangle\dv\\
\nonumber&+\int_M\langle V,R^N(\tau(\phi),\bar\nabla_{e_j}\bar\nabla_{e_k}\big(R^N(V,d\phi(e_k))\tau(\phi)\big))d\phi(e_j)
\rangle\dv\\
\nonumber&+\int_M\langle V,R^N(\tau(\phi),\bar\nabla_{e_j}\bar\Delta^2V)d\phi(e_j)
\rangle\dv\\
\nonumber&-\int_M\langle V,R^N(\tau(\phi),\bar\nabla_{e_j}\bar\Delta\big(R^N(V,d\phi(e_k))d\phi(e_k)\big))d\phi(e_j)
\rangle\dv\\
\nonumber&-\int_M\langle V,R^N(\tau(\phi),\bar\nabla_{e_j}\bar\Delta\tau(\phi))\bar\nabla_{e_j}V
\rangle\dv\\
\nonumber&+\int_M\langle V,(\nabla_{V}R^N)(\bar\nabla_{e_j}\tau(\phi),\bar\Delta\tau(\phi))d\phi(e_j)\rangle\dv
\\
\nonumber&-\int_M\langle V,R^N(\bar\nabla_{e_j}\bar\Delta V,\bar\Delta\tau(\phi))d\phi(e_j)\rangle\dv
\\
\nonumber&+\int_M\langle V,R^N(\bar\nabla_{e_j}\big(R^N(V,d\phi(e_k))d\phi(e_k)\big),\bar\Delta\tau(\phi))d\phi(e_j)\rangle\dv
\\
\nonumber&+\int_M\langle V,R^N(R^N(V,d\phi(e_j))\tau(\phi),\bar\Delta\tau(\phi))d\phi(e_j)\rangle\dv
\\
\nonumber&-\int_M\langle V,R^N(\bar\nabla_{e_j}\tau(\phi),R^N(V,d\phi(e_k))\bar\nabla_{e_k}\tau(\phi))d\phi(e_j)\rangle\dv
\\
\nonumber&-\int_M\langle V,R^N(\bar\nabla_{e_j}\tau(\phi),\bar\nabla_{e_k}\big(R^N(V,d\phi(e_k))\tau(\phi)\big))d\phi(e_j)\rangle\dv
\\
\nonumber&+\int_M\langle V,R^N(\bar\nabla_{e_j}\tau(\phi),\bar\Delta\big(R^N(V,d\phi(e_k))d\phi(e_k)\big))d\phi(e_j)\rangle\dv
\\
\nonumber&-\int_M\langle V,R^N(\bar\nabla_{e_j}\tau(\phi),\bar\Delta^2V)d\phi(e_j)\rangle\dv
\\
\nonumber&+\int_M\langle V,R^N(\bar\nabla_{e_j}\tau(\phi),\bar\Delta\tau(\phi))\bar\nabla_{e_j}V\rangle\dv.
\end{align}

\end{Lem}
\begin{proof}
We consider a variation of the map \(\phi\) as detailed in \eqref{dfn:variation-phi}.
Differentiating \eqref{first-variation-4-energy} once more with respect to \(t\) we obtain
\begin{align}
\label{eq:sv4a}
\frac{d^2}{dt^2}E_4(\phi_t)=&-\int_M\langle\frac{\bar\nabla}{\partial t} d\phi_t(\partial_t),\underbrace{\tau_4(\phi_t)}_{=0}\rangle\dv
-\int_M\langle d\phi_t(\partial_t),\frac{\bar\nabla}{\partial t}\bar\Delta^{3}\tau(\phi_t)\rangle\dv \\
\nonumber&+\int_M\langle d\phi_t(\partial_t),\frac{\bar\nabla}{\partial t}\big(R^N(\bar\Delta^{2}\tau(\phi_t),d\phi_t(e_j))d\phi_t(e_j)\big)\rangle\dv \\
\nonumber&-\int_M\langle d\phi_t(\partial_t),\frac{\bar\nabla}{\partial t}\big(R^N(\tau(\phi_t),\bar\nabla_{e_j}\bar\Delta\tau(\phi_t))d\phi_t(e_j)\big)\rangle\dv \\
\nonumber&+\int_M\langle d\phi_t(\partial_t),\frac{\bar\nabla}{\partial t}\big(R^N(\bar\nabla_{e_j}\tau(\phi_t),\bar\Delta\tau(\phi_t))d\phi_t(e_j)\big)\rangle\dv.
\end{align}
Note that the first term vanishes as we are considering the second variation of a \(4\)-harmonic map. In order to manipulate the second term of \eqref{eq:sv4a} we employ the identities given in
\eqref{commutator-second-variation} and find
\begin{align*}
\frac{\bar\nabla}{\partial t}\bar\Delta^{3}\tau(\phi_t)
=&-R^N(d\phi_t(\partial_t),d\phi_t(e_j))\bar\nabla_{e_j}\bar\Delta^{2}\tau(\phi_t)
-\bar\nabla_{e_j}\big(R^N(d\phi_t(\partial_t),d\phi_t(e_j))\bar\Delta^{2}\tau(\phi_t)\big) \\
&-\bar\Delta\big(R^N(d\phi_t(\partial_t),d\phi_t(e_j))\bar\nabla_{e_j}\bar\Delta\tau(\phi_t)\big)
-\bar\Delta\bar\nabla_{e_j}\big(R^N(d\phi_t(\partial_t),d\phi_t(e_j))\bar\Delta\tau(\phi_t)\big) \\
&-\bar\Delta^2\big(R^N(d\phi_t(\partial_t),d\phi_t(e_j))\bar\nabla_{e_j}\tau(\phi_t)\big)
-\bar\Delta^2\bar\nabla_{e_j}\big(R^N(d\phi_t(\partial_t),d\phi_t(e_j))\tau(\phi_t)\big) \\
&-\bar\Delta^4d\phi_t(\partial_t)
+\bar\Delta^3\big(R^N(d\phi_t(\partial_t),d\phi_t(e_k))d\phi_t(e_k)\big).
\end{align*}
Hence, we obtain
\begin{align*}
\int_M\langle d\phi_t(\partial_t),\frac{\bar\nabla}{\partial t}\bar\Delta^{3}\tau(\phi_t)\rangle\dv\big|_{t=0}
=&-\int_M\langle V,R^N(V,d\phi(e_j))\bar\nabla_{e_j}\bar\Delta^{2}\tau(\phi)
\rangle\dv\\
&+\int_M\langle \bar\nabla_{e_j} V,R^N(V,d\phi(e_j))\bar\Delta^{2}\tau(\phi) 
\rangle\dv\\
&-\int_M\langle \bar\Delta V,R^N(V,d\phi(e_j))\bar\nabla_{e_j}\bar\Delta\tau(\phi)
\rangle\dv\\
&+\int_M\langle \bar\nabla_{e_j}\bar\Delta V,R^N(V,d\phi(e_j))\bar\Delta\tau(\phi)
\rangle\dv\\
&-\int_M\langle \bar\Delta^2 V,R^N(V,d\phi(e_j))\bar\nabla_{e_j}\tau(\phi)
\rangle\dv\\
&+\int_M\langle \bar\nabla_{e_j}\bar\Delta^2 V,R^N(V,d\phi(e_j))\tau(\phi)
\rangle\dv\\
&-\int_M|\bar\Delta^2V|^2\dv\\
&+\int_M\langle \bar\Delta^3 V,R^N(V,d\phi(e_k))d\phi(e_k)
\rangle\dv.
\end{align*}

Concerning the third term of \eqref{eq:sv4a} we find
\begin{align*}
\frac{\bar\nabla}{\partial t}\big(R^N(\bar\Delta^{2}\tau(\phi_t),d\phi_t(e_j))d\phi_t(e_j)\big)=
&(\nabla_{d\phi_t(\partial_t)}R^N)(\bar\Delta^{2}\tau(\phi_t),d\phi_t(e_j))d\phi_t(e_j)\\
&+R^N(\frac{\bar\nabla}{\partial t}\bar\Delta^{2}\tau(\phi_t),d\phi_t(e_j))d\phi_t(e_j)\\
&+R^N(\bar\Delta^{2}\tau(\phi_t),\bar\nabla_{e_j}d\phi_t(\partial_t))d\phi_t(e_j)\\
&+R^N(\bar\Delta^{2}\tau(\phi_t),d\phi_t(e_j))\bar\nabla_{e_j}d\phi_t(\partial_t)\\
=
&(\nabla_{d\phi_t(\partial_t)}R^N)(\bar\Delta^{2}\tau(\phi_t),d\phi_t(e_j))d\phi_t(e_j)\\
&-R^N(R^N(d\phi_t(\partial_t),d\phi_t(e_k))\bar\nabla_{e_k}\bar\Delta\tau(\phi_t),d\phi_t(e_j))d\phi_t(e_j)\\
&-R^N(\bar\nabla_{e_k}\big(R^N(d\phi_t(\partial_t),d\phi_t(e_k))\bar\Delta\tau(\phi_t)\big),d\phi_t(e_j))d\phi_t(e_j)\\
&-R^N(\bar\Delta\big(R^N(d\phi_t(\partial_t),d\phi_t(e_k))\bar\nabla_{e_k}\tau(\phi_t)\big),d\phi_t(e_j))d\phi_t(e_j)\\
&-R^N(\bar\Delta\bar\nabla_{e_k}\big(R^N(d\phi_t(\partial_t),d\phi_t(e_k))\tau(\phi_t)\big),d\phi_t(e_j))d\phi_t(e_j)\\
&-R^N(\bar\Delta^3d\phi_t(\partial_t),d\phi_t(e_j))d\phi_t(e_j)\\
&+R^N(\bar\Delta^{2}\big(R^N(d\phi_t(\partial_t),d\phi_t(e_k))d\phi_t(e_k)\big),d\phi_t(e_j))d\phi_t(e_j)\\
&+R^N(\bar\Delta^{2}\tau(\phi_t),\bar\nabla_{e_j}d\phi_t(\partial_t))d\phi_t(e_j)\\
&+R^N(\bar\Delta^{2}\tau(\phi_t),d\phi_t(e_j))\bar\nabla_{e_j}d\phi_t(\partial_t),
\end{align*}
where we used the identities given by \eqref{commutator-second-variation}.
Hence, we obtain
\begin{align*}
\int_M\langle& d\phi_t(\partial_t),\frac{\bar\nabla}{\partial t}\big(R^N(\bar\Delta^{2}\tau(\phi_t),d\phi_t(e_j))d\phi_t(e_j)\big)\rangle\dv\big|_{t=0} \\
=&\int_M\langle V,(\nabla_{V}R^N)(\bar\Delta^{2}\tau(\phi),d\phi(e_j))d\phi(e_j)
\rangle\dv \\
&-\int_M\langle V,R^N(R^N(V,d\phi(e_k))\bar\nabla_{e_k}\bar\Delta\tau(\phi),d\phi(e_j))d\phi(e_j)
\rangle\dv \\
&-\int_M\langle V,R^N(\bar\nabla_{e_k}\big(R^N(V,d\phi(e_k))\bar\Delta\tau(\phi)\big),d\phi(e_j))d\phi(e_j)
\rangle\dv \\
&-\int_M\langle V,R^N(\bar\Delta\big(R^N(V,d\phi(e_k))\bar\nabla_{e_k}\tau(\phi)\big),d\phi(e_j))d\phi(e_j)
\rangle\dv \\
&-\int_M\langle V,R^N(\bar\Delta\bar\nabla_{e_k}\big(R^N(V,d\phi(e_k))\tau(\phi)\big),d\phi(e_j))d\phi(e_j)
\rangle\dv \\
&-\int_M\langle V,R^N(\bar\Delta^3V,d\phi(e_j))d\phi(e_j)
\rangle\dv \\
&+\int_M\langle V,R^N(\bar\Delta^{2}\big(R^N(V,d\phi(e_k))d\phi(e_k)\big),d\phi(e_j))d\phi(e_j)
\rangle\dv \\
&+\int_M\langle V,R^N(\bar\Delta^{2}\tau(\phi),\bar\nabla_{e_j}V)d\phi(e_j)
\rangle\dv\\
&+\int_M\langle V,R^N(\bar\Delta^{2}\tau(\phi),d\phi(e_j))\bar\nabla_{e_j}V
\rangle\dv
\end{align*}
completing the variation of the third term in \eqref{eq:sv4a}.

Regarding the fourth term of \eqref{eq:sv4a} we start by calculating
\begin{align*}
\frac{\bar\nabla}{\partial t}\big(R^N(\tau(\phi_t),\bar\nabla_{e_j}\bar\Delta\tau(\phi_t))d\phi_t(e_j)\big)=&(\nabla_{d\phi_t(\partial_t)}R^N)(\tau(\phi_t),\bar\nabla_{e_j}\bar\Delta\tau(\phi_t))d\phi_t(e_j)\\
&+R^N(\frac{\bar\nabla}{\partial t}\tau(\phi_t),\bar\nabla_{e_j}\bar\Delta\tau(\phi_t))d\phi_t(e_j)
\\
&+R^N(\tau(\phi_t),\frac{\bar\nabla}{\partial t}\bar\nabla_{e_j}\bar\Delta\tau(\phi_t))d\phi_t(e_j)
\\
&+R^N(\tau(\phi_t),\bar\nabla_{e_j}\bar\Delta\tau(\phi_t))\bar\nabla_{e_j}d\phi_t(\partial_t)\\
=&
(\nabla_{d\phi_t(\partial_t)}R^N)(\tau(\phi_t),\bar\nabla_{e_j}\bar\Delta\tau(\phi_t))d\phi_t(e_j)\\
&-R^N(\bar\Delta d\phi_t(\partial_t),\bar\nabla_{e_j}\bar\Delta\tau(\phi_t))d\phi_t(e_j) \\
&+R^N(R^N(d\phi_t(\partial_t),d\phi_t(e_k))d\phi_t(e_k),\bar\nabla_{e_j}\bar\Delta\tau(\phi_t))d\phi_t(e_j) \\
&+R^N(\tau(\phi_t),R^N(d\phi_t(\partial_t),d\phi_t(e_j))\bar\Delta\tau(\phi_t))d\phi_t(e_j)
\\
&-R^N(\tau(\phi_t),\bar\nabla_{e_j}\big(R^N(d\phi_t(\partial_t),d\phi_t(e_k))\bar\nabla_{e_k}\tau(\phi_t)\big))d\phi_t(e_j)
\\
&-R^N(\tau(\phi_t),\bar\nabla_{e_j}\bar\nabla_{e_k}\big(R^N(d\phi_t(\partial_t),d\phi_t(e_k))\tau(\phi_t)\big))d\phi_t(e_j)
\\
&-R^N(\tau(\phi_t),\bar\nabla_{e_j}\bar\Delta^2d\phi_t(\partial_t))d\phi_t(e_j)
\\
&+R^N(\tau(\phi_t),\bar\nabla_{e_j}\bar\Delta\big(R^N(d\phi_t(\partial_t),d\phi_t(e_k))d\phi_t(e_k)\big))d\phi_t(e_j)
\\
&+R^N(\tau(\phi_t),\bar\nabla_{e_j}\bar\Delta\tau(\phi_t))\bar\nabla_{e_j}d\phi_t(\partial_t).
\end{align*}
Hence, we find, after evaluating at \(t=0\), that
\begin{align*}
\int_M\langle& d\phi_t(\partial_t),\frac{\bar\nabla}{\partial t}\big(R^N(\tau(\phi_t),\bar\nabla_{e_j}\bar\Delta\tau(\phi_t))d\phi_t(e_j)\big)\rangle\dv\big|_{t=0}\\
=&
\int_M\langle V,(\nabla_{V}R^N)(\tau(\phi),\bar\nabla_{e_j}\bar\Delta\tau(\phi))d\phi(e_j)
\rangle\dv\\
&-\int_M\langle V,R^N(\bar\Delta V,\bar\nabla_{e_j}\bar\Delta\tau(\phi))d\phi(e_j)
\rangle\dv\\
&+\int_M\langle V,R^N(R^N(V,d\phi(e_k))d\phi(e_k),\bar\nabla_{e_j}\bar\Delta\tau(\phi))d\phi(e_j)
\rangle\dv\\
&+\int_M\langle V,R^N(\tau(\phi),R^N(V,d\phi(e_j))\bar\Delta\tau(\phi))d\phi(e_j)
\rangle\dv\\
&-\int_M\langle V,R^N(\tau(\phi),\bar\nabla_{e_j}\big(R^N(V,d\phi(e_k))\bar\nabla_{e_k}\tau(\phi)\big))d\phi(e_j)
\rangle\dv\\
&-\int_M\langle V,R^N(\tau(\phi),\bar\nabla_{e_j}\bar\nabla_{e_k}\big(R^N(V,d\phi(e_k))\tau(\phi)\big))d\phi(e_j)
\rangle\dv\\
&-\int_M\langle V,R^N(\tau(\phi),\bar\nabla_{e_j}\bar\Delta^2V)d\phi(e_j)
\rangle\dv\\
&+\int_M\langle V,R^N(\tau(\phi),\bar\nabla_{e_j}\bar\Delta\big(R^N(V,d\phi(e_k))d\phi(e_k)\big))d\phi(e_j)
\rangle\dv\\
&+\int_M\langle V,R^N(\tau(\phi),\bar\nabla_{e_j}\bar\Delta\tau(\phi))\bar\nabla_{e_j}V
\rangle\dv.
\end{align*}
Finally, regarding the last term of \eqref{eq:sv4a} we first calculate
\begin{align*}
\frac{\bar\nabla}{\partial t}\big(R^N(\bar\nabla_{e_j}\tau(\phi_t),\bar\Delta\tau(\phi_t))d\phi_t(e_j)\big)=
&(\nabla_{d\phi_t(\partial_t)}R^N)(\bar\nabla_{e_j}\tau(\phi_t),\bar\Delta\tau(\phi_t))d\phi_t(e_j) \\
&+R^N(\frac{\bar\nabla}{\partial t}\bar\nabla_{e_j}\tau(\phi_t),\bar\Delta\tau(\phi_t))d\phi_t(e_j)
\\
&+R^N(\bar\nabla_{e_j}\tau(\phi_t),\frac{\bar\nabla}{\partial t}\bar\Delta\tau(\phi_t))d\phi_t(e_j)
\\
&+R^N(\bar\nabla_{e_j}\tau(\phi_t),\bar\Delta\tau(\phi_t))\bar\nabla_{e_j}d\phi_t(\partial_t)
\\
=&
(\nabla_{d\phi_t(\partial_t)}R^N)(\bar\nabla_{e_j}\tau(\phi_t),\bar\Delta\tau(\phi_t))d\phi_t(e_j)
\\
&-R^N(\bar\nabla_{e_j}\bar\Delta d\phi_t(\partial_t),\bar\Delta\tau(\phi_t))d\phi_t(e_j)\\
&+R^N(\bar\nabla_{e_j}\big(R^N(d\phi_t(\partial_t),d\phi_t(e_k))d\phi_t(e_k)\big),\bar\Delta\tau(\phi_t))d\phi_t(e_j)\\
&+R^N(R^N(d\phi_t(\partial_t),d\phi_t(e_j))\tau(\phi_t),\bar\Delta\tau(\phi_t))d\phi_t(e_j)\\
&-R^N(\bar\nabla_{e_j}\tau(\phi_t),R^N(d\phi_t(\partial_t),d\phi_t(e_k))\bar\nabla_{e_k}\tau(\phi_t))d\phi_t(e_j)
\\
&-R^N(\bar\nabla_{e_j}\tau(\phi_t),\bar\nabla_{e_k}\big(R^N(d\phi_t(\partial_t),d\phi_t(e_k))\tau(\phi_t)\big))d\phi_t(e_j)
\\
&+R^N(\bar\nabla_{e_j}\tau(\phi_t),\bar\Delta\big(R^N(d\phi_t(\partial_t),d\phi_t(e_k))d\phi_t(e_k)\big))d\phi_t(e_j)
\\
&-R^N(\bar\nabla_{e_j}\tau(\phi_t),\bar\Delta^2d\phi_t(\partial_t))d\phi_t(e_j)
\\
&+R^N(\bar\nabla_{e_j}\tau(\phi_t),\bar\Delta\tau(\phi_t))\bar\nabla_{e_j}d\phi_t(\partial_t).
\end{align*}
Hence, we obtain 
\begin{align*}
\int_M&\langle d\phi_t(\partial_t),\frac{\bar\nabla}{\partial t}\big(R^N(\bar\nabla_{e_j}\tau(\phi_t),\bar\Delta\tau(\phi_t))d\phi_t(e_j)\big)\rangle\dv\big|_{t=0} \\
=&\int_M\langle V,(\nabla_{V}R^N)(\bar\nabla_{e_j}\tau(\phi),\bar\Delta\tau(\phi))d\phi(e_j)\rangle\dv
\\
&-\int_M\langle V,R^N(\bar\nabla_{e_j}\bar\Delta V,\bar\Delta\tau(\phi))d\phi(e_j)\rangle\dv
\\
&+\int_M\langle V,R^N(\bar\nabla_{e_j}\big(R^N(V,d\phi(e_k))d\phi(e_k)\big),\bar\Delta\tau(\phi))d\phi(e_j)\rangle\dv
\\
&+\int_M\langle V,R^N(R^N(V,d\phi(e_j))\tau(\phi),\bar\Delta\tau(\phi))d\phi(e_j)\rangle\dv
\\
&-\int_M\langle V,R^N(\bar\nabla_{e_j}\tau(\phi),R^N(V,d\phi(e_k))\bar\nabla_{e_k}\tau(\phi))d\phi(e_j)\rangle\dv
\\
&-\int_M\langle V,R^N(\bar\nabla_{e_j}\tau(\phi),\bar\nabla_{e_k}\big(R^N(V,d\phi(e_k))\tau(\phi)\big))d\phi(e_j)\rangle\dv
\\
&+\int_M\langle V,R^N(\bar\nabla_{e_j}\tau(\phi),\bar\Delta\big(R^N(V,d\phi(e_k))d\phi(e_k)\big))d\phi(e_j)\rangle\dv
\\
&-\int_M\langle V,R^N(\bar\nabla_{e_j}\tau(\phi),\bar\Delta^2V)d\phi(e_j)\rangle\dv
\\
&+\int_M\langle V,R^N(\bar\nabla_{e_j}\tau(\phi),\bar\Delta\tau(\phi))\bar\nabla_{e_j}V\rangle\dv.
\end{align*}
The claim now follows from adding up all terms.
\end{proof}

Using the symmetries of the Riemann curvature tensor and also employing integration by parts
we are able to simplify the previous formula for the second variation of the 4-energy
significantly.

\begin{Prop}
[Second variation of the 4-energy (simplified version)]
Let \(\phi\colon M\to N\) be a smooth \(4\)-harmonic map.
Consider a variation of the map \(\phi\) as described in \eqref{dfn:variation-phi}.
Then, the second variation of the 4-energy can be simplified to
\begin{align}
\label{sv-4energy-simplified}
\frac{d^2}{dt^2}\big|_{t=0}E_4(\phi_t)
=&
\int_M\big|\bar\Delta^2V-\bar\Delta\big(R^N(V,d\phi(e_k))d\phi(e_k)\big)\big|^2\dv 
\\&\nonumber
+2\int_M\langle R^N(\bar\Delta V,\bar\nabla_{e_j}\bar\Delta\tau(\phi))d\phi(e_j),V\rangle\dv
\\&\nonumber
-2\int_M\langle R^N(\bar\nabla_{e_j}\bar\Delta V,\bar\Delta\tau(\phi))d\phi(e_j),V\rangle\dv
\\&\nonumber
+2\int_M\langle R^N(V,d\phi(e_j))\bar\nabla_{e_j}\tau(\phi),\bar\Delta^2 V\rangle\dv
\\&\nonumber
+2\int_M\langle R^N(\tau(\phi),\bar\nabla_{e_j}\bar\Delta^2V)d\phi(e_j),V
\rangle\dv
\\&\nonumber
-2\int_M\langle R^N(V,d\phi(e_k))d\phi(e_k),R^N(V,d\phi(e_j))\bar\nabla_{e_j}\bar\Delta\tau(\phi)
\rangle\dv
\\&\nonumber
+2\int_M\langle R^N(V,d\phi(e_j))\tau(\phi),R^N(V,d\phi(e_j))\bar\Delta\tau(\phi)
\rangle\dv
\\&\nonumber
-2\int_M\langle R^N(V,d\phi(e_j))d\phi(e_j),\bar\nabla_{e_k}\big(R^N(V,d\phi(e_k))\bar\Delta\tau(\phi)\big)
\rangle\dv
\\&\nonumber
-2\int_M\langle \bar\Delta\big(R^N(V,d\phi(e_j))d\phi(e_j)\big),R^N(V,d\phi(e_k))\bar\nabla_{e_k}\tau(\phi)\rangle\dv
\\&\nonumber
+2\int_M\langle R^N(V,d\phi(e_j))\tau(\phi),\bar\nabla_{e_j}\bar\Delta\big(R^N(V,d\phi(e_k))d\phi(e_k)\big)
\rangle\dv
\\&\nonumber
+2\int_M\langle R^N(d\phi(e_j),V)\tau(\phi),\bar\nabla_{e_j}\big(R^N(V,d\phi(e_k))\bar\nabla_{e_k}\tau(\phi)\big)\rangle\dv \\
\nonumber&-2\int_M\langle R^N(V,d\phi(e_j))\bar\Delta^{2}\tau(\phi),\bar\nabla_{e_j} V
\rangle\dv\\
&+\nonumber
\int_M\big|\bar\nabla_{e_k}\big(R^N(V,d\phi(e_k))\tau(\phi)\big)\big|^2\dv
\\&\nonumber
+\int_M|R^N(V,d\phi(e_k))\bar\nabla_{e_k}\tau(\phi)|^2
\dv\\
\nonumber&+\int_M\langle R^N(V,d\phi(e_j))\bar\nabla_{e_j}\bar\Delta^{2}\tau(\phi),V
\rangle\dv\\
\nonumber&+\int_M\langle R^N(\bar\Delta^{2}\tau(\phi),d\phi(e_j))\bar\nabla_{e_j} V,V
\rangle\dv\\
\nonumber&-\int_M\langle R^N(\tau(\phi),\bar\nabla_{e_j}\bar\Delta\tau(\phi))\bar\nabla_{e_j}V,V
\rangle\dv\\
\nonumber&+\int_M\langle R^N(\bar\nabla_{e_j}\tau(\phi),\bar\Delta\tau(\phi))\bar\nabla_{e_j}V,V
\rangle\dv\\
\nonumber&+\int_M\langle (\nabla_{V}R^N)(\bar\Delta^{2}\tau(\phi),d\phi(e_j))d\phi(e_j),V
\rangle\dv \\
\nonumber&-\int_M\langle (\nabla_{V}R^N)(\tau(\phi),\bar\nabla_{e_j}\bar\Delta\tau(\phi))d\phi(e_j)
,V\rangle\dv\\
\nonumber&+\int_M\langle (\nabla_{V}R^N)(\bar\nabla_{e_j}\tau(\phi),\bar\Delta\tau(\phi))d\phi(e_j),V\rangle\dv.
\end{align}

\end{Prop}

\begin{proof}
First of all, we note that using the symmetries of the Riemann curvature tensor
and due to the self-adjointness of the Laplacian we have
\begin{align*}
\int_M\langle V,R^N(\bar\Delta^{2}\big(R^N(V,d\phi(e_k))d\phi(e_k)\big),d\phi(e_j))d\phi(e_j)
\rangle\dv=\int_M\big|\bar\Delta\big(R^N(V,d\phi(e_k))d\phi(e_k)\big)\big|^2\dv.
\end{align*}
Hence, we can complete the square as follows
\begin{align*}
\int_M&|\bar\Delta^2V|^2\dv-\int_M\langle \bar\Delta^3 V,R^N(V,d\phi(e_k))d\phi(e_k)
\rangle\dv
-\int_M\langle V,R^N(\bar\Delta^3V,d\phi(e_j))d\phi(e_j)
\rangle\dv \\
&+\int_M\langle V,R^N(\bar\Delta^{2}\big(R^N(V,d\phi(e_k))d\phi(e_k)\big),d\phi(e_j))d\phi(e_j)\rangle\dv \\
&=\int_M\big|\bar\Delta^2V-\bar\Delta\big(R^N(V,d\phi(e_k))d\phi(e_k)\big)\big|^2\dv
\end{align*}
giving rise to the first term.

Now, we note that
\begin{align*}
\int_M\langle& \bar\Delta V,R^N(V,d\phi(e_j))\bar\nabla_{e_j}\bar\Delta\tau(\phi)
\rangle\dv
+\int_M\langle V,R^N(\bar\Delta V,\bar\nabla_{e_j}\bar\Delta\tau(\phi))d\phi(e_j)
\rangle\dv\\
&=2\int_M\langle R^N(\bar\Delta V,\bar\nabla_{e_j}\bar\Delta\tau(\phi))d\phi(e_j),V\rangle\dv
\end{align*}
due to the symmetries of the Riemann curvature tensor.
By the same reasoning
\begin{align*}
-\int_M\langle& \bar\nabla_{e_j}\bar\Delta V,R^N(V,d\phi(e_j))\bar\Delta\tau(\phi)\big)
\rangle\dv
-\int_M\langle V,R^N(\bar\nabla_{e_j}\bar\Delta V,\bar\Delta\tau(\phi))d\phi(e_j)\rangle\dv
\\
&=-2\int_M\langle R^N(\bar\nabla_{e_j}\bar\Delta V,\bar\Delta\tau(\phi))d\phi(e_j),V\rangle\dv
\end{align*}
and also
\begin{align*}
\int_M\langle& \bar\Delta^2 V,R^N(V,d\phi(e_j))\bar\nabla_{e_j}\tau(\phi)
\rangle\dv
-\int_M\langle V,R^N(\bar\nabla_{e_j}\tau(\phi),\bar\Delta^2V)d\phi(e_j)\rangle\dv \\
&=2\int_M\langle R^N(V,d\phi(e_j))\bar\nabla_{e_j}\tau(\phi),\bar\Delta^2 V
\rangle\dv.
\end{align*}
Moreover, we find
\begin{align*}
-\int_M\langle & \bar\nabla_{e_j}\bar\Delta^2 V,R^N(V,d\phi(e_j))\tau(\phi)\rangle\dv
+\int_M\langle V,R^N(\tau(\phi),\bar\nabla_{e_j}\bar\Delta^2V)d\phi(e_j)
\rangle\dv\\
&=2\int_M\langle R^N(\tau(\phi),\bar\nabla_{e_j}\bar\Delta^2V)d\phi(e_j),V
\rangle\dv.
\end{align*}
In addition, we get
\begin{align*}
&-\int_M\langle V,R^N(R^N(V,d\phi(e_k))\bar\nabla_{e_k}\bar\Delta\tau(\phi),d\phi(e_j))d\phi(e_j)
\rangle\dv\\
&-\int_M\langle V,R^N(R^N(V,d\phi(e_k))d\phi(e_k),\bar\nabla_{e_j}\bar\Delta\tau(\phi))d\phi(e_j)
\rangle\dv\\
=&-2\int_M\langle R^N(V,d\phi(e_k))d\phi(e_k),R^N(V,d\phi(e_j))\bar\nabla_{e_j}\bar\Delta\tau(\phi)\rangle\dv
\end{align*}
and also
\begin{align*}
\nonumber&-\int_M\langle V,R^N(\bar\nabla_{e_k}\big(R^N(V,d\phi(e_k))\bar\Delta\tau(\phi)\big),d\phi(e_j))d\phi(e_j)
\rangle\dv \\
\nonumber&+\int_M\langle V,R^N(\bar\nabla_{e_j}\big(R^N(V,d\phi(e_k))d\phi(e_k)\big),\bar\Delta\tau(\phi))d\phi(e_j)\rangle\dv \\
&=-2\int_M\langle R^N(V,d\phi(e_j))d\phi(e_j),\bar\nabla_{e_k}\big(R^N(V,d\phi(e_k))\bar\Delta\tau(\phi)\big)\rangle\dv.
\end{align*}

Similarly, we obtain
\begin{align*}
&-\int_M\langle V,R^N(\tau(\phi),R^N(V,d\phi(e_j))\bar\Delta\tau(\phi))d\phi(e_j)
\rangle\dv\\
&+\int_M\langle V,R^N(R^N(V,d\phi(e_j))\tau(\phi),\bar\Delta\tau(\phi))d\phi(e_j)\rangle\dv
\\
=&2\int_M\langle R^N(V,d\phi(e_j))\tau(\phi),R^N(V,d\phi(e_j))\bar\Delta\tau(\phi)\rangle\dv.
\end{align*}
Using the self-adjointness of the Laplacian and the symmetries of the Riemann curvature tensor
we get
\begin{align*}
&-\int_M\langle V,R^N(\bar\Delta\big(R^N(V,d\phi(e_k))\bar\nabla_{e_k}\tau(\phi)\big),d\phi(e_j))d\phi(e_j)
\rangle\dv \\
&+\int_M\langle V,R^N(\bar\nabla_{e_j}\tau(\phi),\bar\Delta\big(R^N(V,d\phi(e_k))d\phi(e_k)\big))d\phi(e_j)\rangle\dv
\\
&=-2\int_M\langle \bar\Delta\big(R^N(V,d\phi(e_j))d\phi(e_j)\big),R^N(V,d\phi(e_k))\bar\nabla_{e_k}\tau(\phi)\rangle\dv.
\end{align*}
Using integration by parts we find
\begin{align*}
&-\int_M\langle V,R^N(\bar\Delta\bar\nabla_{e_k}\big(R^N(V,d\phi(e_k))\tau(\phi)\big),d\phi(e_j))d\phi(e_j)\rangle\dv \\
&-\int_M\langle V,R^N(\tau(\phi),\bar\nabla_{e_j}\bar\Delta\big(R^N(V,d\phi(e_k))d\phi(e_k)\big))d\phi(e_j)
\rangle\dv\\
&=2\int_M\langle R^N(V,d\phi(e_j))\tau(\phi),\bar\nabla_{e_j}\bar\Delta\big(R^N(V,d\phi(e_k))d\phi(e_k)\big)
\rangle\dv.
\end{align*}

Also, we have
\begin{align*}
&\int_M\langle V,R^N(\tau(\phi),\bar\nabla_{e_j}\big(R^N(V,d\phi(e_k))\bar\nabla_{e_k}\tau(\phi)\big))d\phi(e_j)
\rangle\dv\\
&-\int_M\langle V,R^N(\bar\nabla_{e_j}\tau(\phi),\bar\nabla_{e_k}\big(R^N(V,d\phi(e_k))\tau(\phi)\big))d\phi(e_j)\rangle\dv
\\
&=2\int_M\langle R^N(d\phi(e_j),V)\tau(\phi),\bar\nabla_{e_j}\big(R^N(V,d\phi(e_k))\bar\nabla_{e_k}\tau(\phi)\big)\rangle\dv.
\end{align*}

Moreover, note that 
\begin{align*}
-\langle V,R^N(\bar\nabla_{e_j}\tau(\phi),R^N(V,d\phi(e_k))\bar\nabla_{e_k}\tau(\phi))d\phi(e_j)\rangle=|R^N(V,d\phi(e_k))\bar\nabla_{e_k}\tau(\phi)|^2
\end{align*}
and also
\begin{align*}
\langle& V,R^N(\bar\Delta^2\tau(\phi),\bar\nabla_{e_j}V)d\phi(e_j)\rangle
-\langle R^N(V,d\phi(e_j))\bar\Delta^2\tau(\phi),\bar\nabla_{e_j}V\rangle \\
&=-2\langle R^N(V,d\phi(e_j))\bar\Delta^2\tau(\phi),\bar\nabla_{e_j}V\rangle.
\end{align*}

Again, using integration by parts we find
\begin{align*}
\int_M&\langle V,R^N(\tau(\phi),\bar\nabla_{e_j}\bar\nabla_{e_k}\big(R^N(V,d\phi(e_k))\tau(\phi)\big))d\phi(e_j)
\rangle\dv \\
&=\int_M\big|\bar\nabla_{e_k}\big(R^N(V,d\phi(e_k))\tau(\phi)\big)\big|^2\dv.
\end{align*}
The claim now follows from combining all identities and inserting these into \eqref{eq:snd-variation-4-rough}.
\end{proof}

A first consequence of the previous calculations is the following result
which was already given in \cite{MR3007953} in the context
of the second variation of \(k\)-harmonic maps.

\begin{Satz}
A harmonic map \(\phi\colon M\to N\) is a weakly stable 4-harmonic map.
 \end{Satz}
\begin{proof}
Using that \(\phi\) is harmonic in \eqref{sv-4energy-simplified} we find
\begin{align*}
\frac{d^2}{dt^2}\big|_{t=0}E_4(\phi_t)
=&
\int_M\big|\bar\Delta^2V-\bar\Delta\big(R^N(V,d\phi(e_k))d\phi(e_k)\big)\big|^2\dv \geq 0,
\end{align*}
which already completes the proof.
\end{proof}

\subsection{Space form target}
From now on, we will assume that the target manifold is a space form of constant curvature \(K\).
This will help us a lot in simplifying the long expression for the second variation of \(E_4(\phi)\) given by \eqref{sv-4energy-simplified}.

\begin{Prop}
Let \(\phi\colon M\to N\) be a smooth \(4\)-harmonic map where \(N\)
is a space form of constant curvature \(K\).
Consider a variation of the map \(\phi\) as described in \eqref{dfn:variation-phi}.
Then, the second variation of the 4-energy simplifies to
\begin{align}
\label{sv-4energy-space-form}
\frac{d^2}{dt^2}\big|_{t=0}E_4(\phi_t)
=&\int_M\big|\bar\Delta^2V-K\bar\Delta\big(|d\phi|^2V-\langle d\phi(e_k),V\rangle d\phi(e_k)\big)\big|^2\dv 
\\&\nonumber
+K\int_M\bigg(
3\langle V,\bar\Delta^2\tau(\phi)\rangle\langle d\phi,\bar\nabla V\rangle
-2\langle d\phi(e_j),\bar\Delta^2\tau(\phi)\rangle\langle V,\bar\nabla_{e_j}V\rangle
\\ \nonumber&
+2\langle\bar\nabla\bar\Delta\tau(\phi),d\phi\rangle\langle\bar\Delta V,V\rangle
-2\langle\bar\Delta V,d\phi(e_j)\rangle\langle \bar\nabla_{e_j}\bar\Delta\tau(\phi),V\rangle 
\\ \nonumber&
-2\langle\bar\Delta\tau(\phi),d\phi(e_j)\rangle\langle \bar\nabla_{e_j}\bar\Delta V,V\rangle
+2\langle\bar\nabla\bar\Delta V,d\phi\rangle\langle\bar\Delta\tau(\phi),V\rangle
\\ \nonumber&
+2\langle d\phi,\bar\nabla\tau(\phi)\rangle\langle V,\bar\Delta^2V\rangle
-2\langle V,\bar\nabla_{e_j}\tau(\phi)\rangle\langle d\phi(e_j),\bar\Delta^2V\rangle
\\ \nonumber&
+2\langle\bar\nabla\bar\Delta^2V,d\phi\rangle\langle\tau(\phi),V\rangle
-2\langle\tau(\phi),d\phi(e_j)\rangle\langle \bar\nabla_{e_j}\bar\Delta^2V,V\rangle
\\ \nonumber&
+\langle d\phi,\bar\nabla\bar\Delta^{2}\tau(\phi)\rangle|V|^2
-\langle V,\bar\nabla_{e_j}\bar\Delta^{2}\tau(\phi)\rangle\langle V,d\phi(e_j)\rangle 
\\ \nonumber&
-\langle\bar\nabla\bar\Delta\tau(\phi),\bar\nabla V\rangle\langle\tau(\phi),V\rangle
+\langle\tau(\phi),\bar\nabla_{e_j}V\rangle
\langle \bar\nabla_{e_j}\bar\Delta\tau(\phi),V\rangle
\\&\nonumber
+\langle\bar\Delta\tau(\phi),\bar\nabla_{e_j}V\rangle\langle\bar\nabla_{e_j}\tau(\phi),V\rangle
-\langle\bar\nabla\tau(\phi),\bar\nabla V\rangle\langle\bar\Delta\tau(\phi),V\rangle
\\ \nonumber&
-\langle\bar\Delta^2\tau(\phi),\bar\nabla_{e_j}V\rangle\langle V,d\phi(e_j)\rangle
\bigg)\dv  
\\&\nonumber
+K^2\int_M\bigg(
-2|d\phi|^2|V|^2\langle d\phi,\bar\nabla\bar\Delta\tau(\phi)\rangle
+2|d\phi|^2\langle V,d\phi(e_j)\rangle\langle V,\bar\nabla_{e_j}\bar\Delta\tau(\phi)\rangle \\
\nonumber&+2|\langle d\phi,V\rangle|^2\langle d\phi,\bar\nabla\bar\Delta V\rangle
-2\langle d\phi(e_k),V\rangle\langle V,\bar\nabla_{e_j}\bar\Delta\tau(\phi)\rangle
\langle d\phi(e_k),d\phi(e_j)\rangle
\\&\nonumber
+2\langle d\phi(e_j),\tau(\phi)\rangle \langle d\phi(e_j),\bar\Delta\tau(\phi)\rangle |V|^2
-2\langle d\phi(e_j),\tau(\phi)\rangle\langle d\phi(e_j),V\rangle\langle V,\bar\Delta\tau(\phi)\rangle \\
\nonumber&-2\langle V,\tau(\phi)\rangle\langle d\phi(e_j),\bar\Delta\tau(\phi)\rangle\langle V,d\phi(e_j)\rangle
+2|d\phi|^2\langle V,\tau(\phi)\rangle\langle V,\bar\Delta\tau(\phi)\rangle
\\&\nonumber
-2|d\phi|^2\langle V,\bar\nabla_{e_k}\big(\langle d\phi(e_k),\bar\Delta\tau(\phi)\rangle V\big)\rangle
+2|d\phi|^2\langle V,\bar\nabla_{e_k}\big(\langle V,\bar\Delta\tau(\phi)\rangle d\phi(e_k)\big)\rangle \\
\nonumber&+2\langle d\phi(e_j),V\rangle \langle d\phi(e_j),\bar\nabla_{e_k}\big(\langle d\phi(e_k),\bar\Delta\tau(\phi)\rangle V\big)\rangle \\
&\nonumber-2\langle d\phi(e_j),V\rangle\langle d\phi(e_j),\bar\nabla_{e_k}\big(\langle V,\bar\Delta\tau(\phi)\rangle d\phi(e_k)\big)\rangle
\\&\nonumber
-2\langle\bar\Delta(|d\phi|^2V),V\rangle\langle d\phi,\bar\nabla\tau(\phi)\rangle
+2\langle\bar\Delta(|d\phi|^2V),d\phi(e_k)\rangle\langle V,\bar\nabla_{e_k}\tau(\phi)\rangle \\
\nonumber&
+2\langle\bar\Delta(\langle d\phi(e_j),V\rangle d\phi(e_j)),V\rangle\langle d\phi,\bar\nabla\tau(\phi)\rangle \\
\nonumber&
-2\langle\bar\Delta(\langle d\phi(e_j),V\rangle d\phi(e_j)),d\phi(e_k)\rangle\langle V,\bar\nabla_{e_k}\tau(\phi)\rangle\\&\nonumber
+2\langle d\phi(e_j),\tau(\phi)\rangle\langle V,\bar\nabla_{e_j}\bar\Delta (|d\phi|^2V)\rangle \\
\nonumber&
-2\langle d\phi(e_j),\tau(\phi)\rangle \langle V,\bar\nabla_{e_j}\bar\Delta\big(\langle d\phi(e_k),V\rangle d\phi(e_k)\big)\rangle \\
\nonumber&-2\langle V,\tau(\phi)\rangle\langle d\phi,\bar\nabla\bar\Delta (|d\phi|^2V)\rangle \\
\nonumber&
+2\langle V,\tau(\phi)\rangle\langle d\phi,\bar\nabla\bar\Delta\big(\langle d\phi(e_k),V\rangle d\phi(e_k)\big)\rangle
\\&\nonumber
+2\langle V,\tau(\phi)\rangle\langle d\phi,\bar\nabla\big(\langle d\phi,\bar\nabla\tau(\phi)\rangle V\big)\rangle
\\&\nonumber
-2\langle V,\tau(\phi)\rangle\langle d\phi,\bar\nabla\big(\langle V,\bar\nabla_{e_k}\tau(\phi)\rangle d\phi(e_k)\big)\rangle \\
\nonumber&-2\langle d\phi(e_j),\tau(\phi)\rangle\langle V,\bar\nabla_{e_j}\big(\langle d\phi,\bar\nabla\tau(\phi)\rangle V\big)\rangle
\\&\nonumber
+2\langle d\phi(e_j),\tau(\phi)\rangle\langle V,\bar\nabla_{e_j}\big(\langle V,\bar\nabla_{e_k}\tau(\phi)\rangle d\phi(e_k)\big)\rangle
\\&\nonumber
+|\bar\nabla_{e_k}\big(\langle d\phi(e_k),\tau(\phi)\rangle V
-\langle V,\tau(\phi)\rangle d\phi(e_k)\big)|^2
\\&\nonumber
+|\langle d\phi,\bar\nabla\tau(\phi)\rangle|^2|V|^2
+|\langle V,\bar\nabla_{e_k}\tau(\phi)\rangle d\phi(e_k)|^2 \\
\nonumber&-2\langle d\phi,\bar\nabla\tau(\phi)\rangle\langle V,\bar\nabla_{e_k}\tau(\phi)\rangle
\langle V,d\phi(e_k)\rangle\bigg)\dv.
\end{align}
\end{Prop}

\begin{proof}
First of all we note that all terms in \eqref{sv-4energy-simplified} involving the derivative of the Riemann curvature tensor drop out as we are considering a space form of constant curvature \(K\).
Now, exploiting the particular form of the Riemann curvature tensor \eqref{curvature-space-form}
we can simplify
\begin{align*}
\langle R^N(\bar\Delta V,\bar\nabla_{e_j}\bar\Delta\tau(\phi))d\phi(e_j),V\rangle=&
K(\langle\bar\nabla\bar\Delta\tau(\phi),d\phi\rangle\langle\bar\Delta V,V\rangle
-\langle\bar\Delta V,d\phi(e_j)\rangle\langle \bar\nabla_{e_j}\bar\Delta\tau(\phi),V\rangle),
\\
\langle R^N(\bar\nabla_{e_j}\bar\Delta V,\bar\Delta\tau(\phi))d\phi(e_j),V\rangle=&
K(\langle\bar\Delta\tau(\phi),d\phi(e_j)\rangle\langle \bar\nabla_{e_j}\bar\Delta V,V\rangle
-\langle\bar\nabla\bar\Delta V,d\phi\rangle\langle\bar\Delta\tau(\phi),V\rangle),
\\
\langle R^N(V,d\phi(e_j))\bar\nabla_{e_j}\tau(\phi),\bar\Delta^2 V\rangle=&
K(\langle d\phi,\bar\nabla\tau(\phi)\rangle\langle V,\bar\Delta^2V\rangle
-\langle V,\bar\nabla_{e_j}\tau(\phi)\rangle\langle d\phi(e_j),\bar\Delta^2V\rangle),
\\
\langle R^N(\tau(\phi),\bar\nabla_{e_j}\bar\Delta^2V)d\phi(e_j),V\rangle=&
K(\langle\bar\nabla\bar\Delta^2V,d\phi\rangle\langle\tau(\phi),V\rangle
-\langle\tau(\phi),d\phi(e_j)\rangle\langle \bar\nabla_{e_j}\bar\Delta^2V,V\rangle),
\\
\langle R^N(V,d\phi(e_j))\bar\nabla_{e_j}\bar\Delta^{2}\tau(\phi),V
\rangle=&
K(\langle d\phi,\bar\nabla\bar\Delta^{2}\tau(\phi)\rangle|V|^2
-\langle V,\bar\nabla_{e_j}\bar\Delta^{2}\tau(\phi)\rangle\langle V,d\phi(e_j)\rangle),
\\
\langle R^N(V,d\phi(e_j))\bar\Delta^{2}\tau(\phi),\bar\nabla_{e_j}V
\rangle=& 
K(\langle d\phi(e_j),\bar\Delta^2\tau(\phi)\rangle\langle V,\bar\nabla_{e_j}V\rangle
-\langle V,\bar\Delta^2\tau(\phi)\rangle\langle d\phi,\bar\nabla V\rangle),
\\
\langle R^N(\tau(\phi),\bar\nabla_{e_j}\bar\Delta\tau(\phi))\bar\nabla_{e_j}V,V
\rangle=&
K(\langle\bar\nabla\bar\Delta\tau(\phi),\bar\nabla V\rangle\langle\tau(\phi),V\rangle
-\langle\tau(\phi),\bar\nabla_{e_j}V\rangle
\langle \bar\nabla_{e_j}\bar\Delta\tau(\phi),V\rangle),
\\
\langle R^N(\bar\nabla_{e_j}\tau(\phi),\bar\Delta\tau(\phi))\bar\nabla_{e_j}V,V\rangle=&
K(\langle\bar\Delta\tau(\phi),\nabla_{e_j}V\rangle\langle\bar\nabla_{e_j}\tau(\phi),V\rangle
-\langle\bar\nabla\tau(\phi),\bar\nabla V\rangle\langle\bar\Delta\tau(\phi),V\rangle),\\
\langle R^N(\bar\Delta^2\tau(\phi),d\phi(e_j))\bar\nabla_{e_j}V,V\rangle
=&K(\langle d\phi,\bar\nabla V\rangle\langle V,\bar\Delta^2\tau(\phi)\rangle
-\langle\bar\Delta^2\tau(\phi),\bar\nabla_{e_j}V\rangle\langle V,d\phi(e_j)\rangle).
\end{align*}

Regarding the terms that are quadratic in the Riemann curvature tensor we obtain
\begin{align*}
\langle &R^N(V,d\phi(e_k))d\phi(e_k),R^N(V,d\phi(e_j))\bar\nabla_{e_j}\bar\Delta\tau(\phi)
\rangle \\
=&K^2\big(|d\phi|^2|V|^2\langle d\phi,\bar\nabla\bar\Delta\tau(\phi)\rangle
-|d\phi|^2\langle V,d\phi(e_j)\rangle\langle V,\bar\nabla_{e_j}\bar\Delta\tau(\phi)\rangle \\
&-|\langle d\phi,V\rangle|^2\langle d\phi,\bar\nabla\bar\Delta V\rangle
+\langle d\phi(e_k),V\rangle\langle V,\bar\nabla_{e_j}\bar\Delta\tau(\phi)\rangle
\langle d\phi(e_k),d\phi(e_j)\rangle\big), \\
\langle &R^N(V,d\phi(e_j))\tau(\phi),R^N(V,d\phi(e_j))\bar\Delta\tau(\phi)
\rangle \\
=&K^2\big(
\langle d\phi(e_j),\tau(\phi)\rangle \langle d\phi(e_j),\bar\Delta\tau(\phi)\rangle |V|^2
-\langle d\phi(e_j),\tau(\phi)\rangle\langle d\phi(e_j),V\rangle\langle V,\bar\Delta\tau(\phi)\rangle \\
&-\langle V,\tau(\phi)\rangle\langle d\phi(e_j),\bar\Delta\tau(\phi)\rangle\langle V,d\phi(e_j)\rangle+|d\phi|^2\langle V,\tau(\phi)\rangle\langle V,\bar\Delta\tau(\phi)\rangle\big).
\end{align*}
In addition, we find
\begin{align*}
\langle &R^N(V,d\phi(e_j))d\phi(e_j),\bar\nabla_{e_k}\big(R^N(V,d\phi(e_k))\bar\Delta\tau(\phi)\big)
\rangle \\
=&K^2\big(|d\phi|^2\langle V,\bar\nabla_{e_k}\big(\langle d\phi(e_k),\bar\Delta\tau(\phi)\rangle V\big)\rangle
-|d\phi|^2\langle V,\bar\nabla_{e_k}\big(\langle V,\bar\Delta\tau(\phi)\rangle d\phi(e_k)\big)\rangle \\
&-\langle d\phi(e_j),V\rangle \langle d\phi(e_j),\bar\nabla_{e_k}\big(\langle d\phi(e_k),\bar\Delta\tau(\phi)\rangle V\big)\rangle 
+\langle d\phi(e_j),V\rangle\langle d\phi(e_j),\bar\nabla_{e_k}\big(\langle V,\bar\Delta\tau(\phi)\rangle d\phi(e_k)\big)\rangle \big)
\end{align*}
and also
\begin{align*}
\langle &\bar\Delta\big(R^N(V,d\phi(e_j))d\phi(e_j)\big),R^N(V,d\phi(e_k))\bar\nabla_{e_k}\tau(\phi)\rangle \\
=&K^2\big(\langle\bar\Delta(|d\phi|^2V),V\rangle\langle d\phi,\bar\nabla\tau(\phi)\rangle
-\langle\bar\Delta(|d\phi|^2V),d\phi(e_k)\rangle\langle V,\bar\nabla_{e_k}\tau(\phi)\rangle 
\\ &
-\langle\bar\Delta(\langle d\phi(e_j),V\rangle d\phi(e_j)),V\rangle\langle d\phi,\bar\nabla\tau(\phi)\rangle
+\langle\bar\Delta(\langle d\phi(e_j),V\rangle d\phi(e_j)),d\phi(e_k)\rangle\langle V,\bar\nabla_{e_k}\tau(\phi)\rangle
\big).
\end{align*}
Moreover, similar manipulations as before yield
\begin{align*}
\langle & R^N(V,d\phi(e_j))\tau(\phi),\bar\nabla_{e_j}\bar\Delta\big(R^N(V,d\phi(e_k))d\phi(e_k)\big)\rangle \\
=& K^2\big(\langle d\phi(e_j),\tau(\phi)\rangle\langle V,\bar\nabla_{e_j}\bar\Delta (|d\phi|^2V)\rangle
-\langle d\phi(e_j),\tau(\phi)\rangle \langle V,\bar\nabla_{e_j}\bar\Delta\big(\langle d\phi(e_k),V\rangle d\phi(e_k)\big)\rangle \\
&-\langle V,\tau(\phi)\rangle\langle d\phi,\bar\nabla\bar\Delta (|d\phi|^2V)\rangle
+\langle V,\tau(\phi)\rangle\langle d\phi,\bar\nabla\bar\Delta\big(\langle d\phi(e_k),V\rangle d\phi(e_k)\big)\rangle
\big)
\end{align*}
and also
\begin{align*}
\langle & R^N(d\phi(e_j),V)\tau(\phi),\bar\nabla_{e_j}\big(R^N(V,d\phi(e_k))\bar\nabla_{e_k}\tau(\phi)\big)\rangle \\
=& K^2\big(
\langle V,\tau(\phi)\rangle\langle d\phi,
\bar\nabla\big(\langle d\phi,\bar\nabla\tau(\phi)\rangle V\big)\rangle
-\langle V,\tau(\phi)\rangle\langle d\phi,
\bar\nabla\big(\langle V,\bar\nabla_{e_k}\tau(\phi)\rangle d\phi(e_k)\big)\rangle \\
&-\langle d\phi(e_j),\tau(\phi)\rangle\langle V,\bar\nabla_{e_j}\big(\langle d\phi,\bar\nabla\tau(\phi)\rangle V\big)\rangle
+\langle d\phi(e_j),\tau(\phi)\langle V,\bar\nabla_{e_j}\big(\langle V,\bar\nabla_{e_k}\tau(\phi)\rangle d\phi(e_k)\big)\rangle.
\end{align*}

Finally, we obtain
\begin{align*}
\big|\bar\nabla_{e_k}\big(R^N(V,d\phi(e_k))\tau(\phi)\big)\big|^2=& 
K^2|\bar\nabla_{e_k}\big(\langle d\phi(e_k),\tau(\phi)\rangle V
-\langle V,\tau(\phi)\rangle d\phi(e_k)\big)|^2,
\\
|R^N(V,d\phi(e_k))\bar\nabla_{e_k}\tau(\phi)|^2=&
K^2\big(|\langle d\phi,\bar\nabla\tau(\phi)\rangle|^2|V|^2
+|\langle V,\bar\nabla_{e_k}\tau(\phi)\rangle d\phi(e_k)|^2 \\
&-2\langle d\phi,\bar\nabla\tau(\phi)\rangle\langle V,\bar\nabla_{e_k}\tau(\phi)\rangle
\langle V,d\phi(e_k)\rangle\big).
\end{align*}
The claim now follows by adding up all contributions.
\end{proof}

\subsection{Normal stability of 4-harmonic hypersurfaces in space forms}
Now, we will further simplify the formula for the second variation 
\eqref{sv-4energy-simplified}
by assuming that
the map \(\phi\) describes a hypersurface in a Riemannian space form of constant curvature \(K\)
and focus on its stability with respect to normal variations.
\begin{Cor}
Let \(\phi\colon M\to N\) be a smooth \(4\)-harmonic hypersurface where \(N\)
is a space form of constant curvature \(K\).
Consider a variation of the map \(\phi\) as described in \eqref{dfn:variation-phi}.
Suppose that \(V=f\nu\), where \(f\in C_c^\infty(M)\) and \(\nu\) is the unit
normal of the hypersurface.
Then, the second variation of the 4-energy simplifies to
\begin{align}
\label{sv-4energy-space-form-hyper}
\frac{d^2}{dt^2}\big|_{t=0}E_4(\phi_t)
=&\int_M\big|\bar\Delta^2(f\nu)-mK\bar\Delta(f\nu)\big|^2\dv 
\\&\nonumber
+K\int_M\bigg(
3\langle f\nu,\bar\Delta^2\tau(\phi)\rangle\langle d\phi,\bar\nabla (f\nu)\rangle
\\&\nonumber
-2\langle d\phi(e_j),\bar\Delta^2\tau(\phi)\rangle\langle f\nu,\bar\nabla_{e_j}(f\nu)\rangle
\\ \nonumber&
+2\langle\bar\nabla\bar\Delta\tau(\phi),d\phi\rangle\langle\bar\Delta (f\nu),f\nu\rangle
\\&\nonumber
-2\langle\bar\Delta (f\nu),d\phi(e_j)\rangle\langle \bar\nabla_{e_j}\bar\Delta\tau(\phi),f\nu\rangle 
\\ \nonumber&
-2\langle\bar\Delta\tau(\phi),d\phi(e_j)\rangle\langle \bar\nabla_{e_j}\bar\Delta (f\nu),f\nu\rangle
\\&\nonumber
+2\langle\bar\nabla\bar\Delta (f\nu),d\phi\rangle\langle\bar\Delta\tau(\phi),f\nu\rangle
\\ \nonumber&
-2|\tau(\phi)|^2\langle f\nu,\bar\Delta^2(f\nu)\rangle
\\&\nonumber
-2\langle f\nu,\bar\nabla_{e_j}\tau(\phi)\rangle\langle d\phi(e_j),\bar\Delta^2(f\nu)\rangle
\\ \nonumber&
+2\langle\bar\nabla\bar\Delta^2(f\nu),d\phi\rangle\langle\tau(\phi),f\nu\rangle
\\ \nonumber&
+\langle d\phi,\bar\nabla\bar\Delta^{2}\tau(\phi)\rangle f^2
\\ \nonumber&
-\langle\bar\nabla\bar\Delta\tau(\phi),\bar\nabla (f\nu)\rangle\langle\tau(\phi),f\nu\rangle
\\&\nonumber
+\langle\tau(\phi),\bar\nabla_{e_j}(f\nu)\rangle
\langle \bar\nabla_{e_j}\bar\Delta\tau(\phi),f\nu\rangle
\\&\nonumber
+\langle\bar\Delta\tau(\phi),\bar\nabla_{e_j}(f\nu)\rangle\langle\bar\nabla_{e_j}\tau(\phi),f\nu\rangle
\\&\nonumber
-\langle\bar\nabla\tau(\phi),\bar\nabla (f\nu)\rangle\langle\bar\Delta\tau(\phi),f\nu\rangle
\bigg)  \dv
\\&\nonumber
+K^2\int_M\bigg(
-2|d\phi|^2f^2\langle d\phi,\bar\nabla\bar\Delta\tau(\phi)\rangle
 \\
\nonumber&
+2|d\phi|^2\langle f\nu,\tau(\phi)\rangle\langle f\nu,\bar\Delta\tau(\phi)\rangle
\\&\nonumber
-2|d\phi|^2\langle f\nu,\bar\nabla_{e_k}\big(\langle d\phi(e_k),\bar\Delta\tau(\phi)\rangle f\nu\big)\rangle
\\&\nonumber
+2|d\phi|^2\langle f\nu,\bar\nabla_{e_k}\big(\langle f\nu,\bar\Delta\tau(\phi)\rangle d\phi(e_k)\big)\rangle 
\\&\nonumber
+2\langle\bar\Delta(|d\phi|^2f\nu),f\nu\rangle|\tau(\phi)|^2
\\&\nonumber
+2\langle\bar\Delta(|d\phi|^2f\nu),d\phi(e_k)\rangle\langle f\nu,\bar\nabla_{e_k}\tau(\phi)\rangle 
\\
\nonumber&-2\langle f\nu,\tau(\phi)\rangle\langle d\phi,\bar\nabla\bar\Delta (|d\phi|^2f\nu)\rangle
\\&\nonumber
+2\langle f\nu,\tau(\phi)\rangle\langle d\phi,\bar\nabla\big(\langle d\phi,\bar\nabla\tau(\phi)\rangle f\nu\big)\rangle
\\&\nonumber
-2\langle f\nu,\tau(\phi)\rangle\langle d\phi,\bar\nabla\big(\langle f\nu,\bar\nabla_{e_k}\tau(\phi)\rangle d\phi(e_k)\big)\rangle \\
\nonumber&
+|\bar\nabla_{e_k}\big(\langle f\nu,\tau(\phi)\rangle d\phi(e_k)\big)|^2
\\&\nonumber
+|\tau(\phi)|^4f^2
\\&\nonumber
+|\langle f\nu,\bar\nabla_{e_k}\tau(\phi)\rangle d\phi(e_k)|^2
\bigg)\dv.
\end{align}
\end{Cor}
\begin{proof}
By assumption \(\phi\) is a hypersurface with unit normal \(\nu\) and thus 
\(\tau(\phi)=mH\nu\) such that
\begin{align*}
\langle d\phi(X),\tau(\phi)\rangle=0
\end{align*}
for all \(X\in TM\). Moreover, we have chosen \(V=f\nu\) and thus also
\begin{align*}
\langle d\phi(X),V\rangle=0
\end{align*}
for all \(X\in TM\). The result is now a direct consequence of equation \eqref{sv-4energy-space-form}.
\end{proof}

Using the formula \eqref{eq:laplace-hypersurface} applied in the case of \(V=f\nu\)
we are then led to the following statement:

\begin{Prop}
Let \(\phi\colon M\to N\) be a smooth \(4\)-harmonic hypersurface with parallel shape operator, i.e. \(\nabla A=0\), where \(N\) is a space form of constant curvature \(K\).
Then, the quadratic form characterizing its normal stability is given by
\begin{align}
\label{qf4-a}
Q_4(f\nu,f\nu)=&\int_M\big|\bar\Delta^2(f\nu)-mK\bar\Delta(f\nu)\big|^2\dv \\
\nonumber&+K\int_M\bigg(-4m^2H^2|\bar\Delta(f\nu)|^2 \\
\nonumber&-4m^2H^2|A|^2|\nabla f|^2 
-10m^2H^2|A|^4f^2
\\
\nonumber&
-4mH\langle A(\operatorname{grad}\Delta f),\operatorname{grad}f\rangle 
-4mH\langle A(e_i),A(\operatorname{grad} f)\rangle\langle A(e_i),\operatorname{grad} f\rangle
\\
\nonumber&-8mH|A|^2\langle A(\operatorname{grad} f),\operatorname{grad} f\rangle 
-4mH\langle A(\Delta\operatorname{grad} f),\operatorname{grad} f\rangle 
\nonumber \bigg)\dv
\\
\nonumber&+K^2\int_M\bigg(
(4m^3H^2+m^2H^2)|\nabla f|^2+(10m^3|A|^2H^2+4m^4H^4)f^2 \\
\nonumber&+ 4m^2H\langle A(\operatorname{grad} f),\operatorname{grad} f\rangle
\bigg)\dv,
\end{align}
where \(f\in C^\infty(M)\) and \(\nu\) represents the normal of the hypersurface.
\end{Prop}

\begin{proof}
We start by manipulating the terms in \eqref{sv-4energy-space-form-hyper} that are 
proportional to \(K\) and find
\begin{align*}
\langle\bar\nabla\bar\Delta\tau(\phi),d\phi\rangle\langle\bar\Delta (f\nu),f\nu\rangle
&=-m^2H^2|A|^2(\Delta ff+|A|^2f^2)
,\\
\langle\bar\Delta (f\nu),d\phi(e_j)\rangle\langle \bar\nabla_{e_j}\bar\Delta\tau(\phi),f\nu\rangle 
&=0
,\\
\langle\bar\Delta\tau(\phi),d\phi(e_j)\rangle\langle \bar\nabla_{e_j}\bar\Delta (f\nu),f\nu\rangle
&=0
,\\
\langle\bar\nabla\bar\Delta (f\nu),d\phi\rangle\langle\bar\Delta\tau(\phi),f\nu\rangle
&=mH|A|^2\langle\bar\nabla\bar\Delta (f\nu),d\phi\rangle f
,\\
|\tau(\phi)|^2\langle f\nu,\bar\Delta^2(f\nu)\rangle
&=m^2H^2\langle f\nu,\bar\Delta^2(f\nu)\rangle
,\\
\langle f\nu,\bar\nabla_{e_j}\tau(\phi)\rangle\langle d\phi(e_j),\bar\Delta^2(f\nu)\rangle
&=0
,\\
\langle\bar\nabla\bar\Delta^2(f\nu),d\phi\rangle\langle\tau(\phi),f\nu\rangle
&=mH\langle\bar\nabla\bar\Delta^2(f\nu),d\phi\rangle f
,\\
\langle d\phi,\bar\nabla\bar\Delta^{2}\tau(\phi)\rangle f^2
&=-m^2H^2|A|^4f^2
,\\
\langle d\phi(e_j),\bar\Delta^2\tau(\phi)\rangle\langle f\nu,\bar\nabla_{e_j}(f\nu)\rangle
&=0
,\\
\langle f\nu,\bar\Delta^2\tau(\phi)\rangle\langle d\phi,\bar\nabla (f\nu)\rangle
&=-m^2H^2|A|^4f^2
,\\
\langle\bar\nabla\bar\Delta\tau(\phi),\bar\nabla (f\nu)\rangle\langle\tau(\phi),f\nu\rangle
&=m^2H^2|A|^4f^2
,\\
\langle\tau(\phi),\bar\nabla_{e_j}(f\nu)\rangle
\langle \bar\nabla_{e_j}\bar\Delta\tau(\phi),f\nu\rangle
&=0
,\\
\langle\bar\Delta\tau(\phi),\bar\nabla_{e_j}(f\nu)\rangle\langle\bar\nabla_{e_j}\tau(\phi),f\nu\rangle
&=0,\\
\langle\bar\nabla\tau(\phi),\bar\nabla (f\nu)\rangle\langle\bar\Delta\tau(\phi),f\nu\rangle
&=m^2H^2|A|^4f^2.
\end{align*}

Concerning the terms that are proportional to \(K^2\) in \eqref{sv-4energy-space-form-hyper}
a direct calculation using \eqref{eq:laplace-hypersurface} yields
\begin{align*}
|d\phi|^2f^2\langle d\phi,\bar\nabla\bar\Delta\tau(\phi)\rangle=&-m^3|A|^2H^2f^2,
 \\
|d\phi|^2\langle f\nu,\tau(\phi)\rangle\langle f\nu,\bar\Delta\tau(\phi)\rangle
=&m^3|A|^2H^2f^2,
\\
|d\phi|^2\langle f\nu,\bar\nabla_{e_k}\big(\langle d\phi(e_k),\bar\Delta\tau(\phi)\rangle f\nu\big)\rangle
=&0,\\
|d\phi|^2\langle f\nu,\bar\nabla_{e_k}\big(\langle f\nu,\bar\Delta\tau(\phi)\rangle d\phi(e_k)\big)\rangle 
=&m^3|A|^2H^2f^2,
\\
\langle\bar\Delta(|d\phi|^2f\nu),f\nu\rangle|\tau(\phi)|^2
=&m^3H^2(\Delta f+|A|^2f)f,\\
\langle\bar\Delta(|d\phi|^2f\nu),d\phi(e_k)\rangle\langle f\nu,\bar\nabla_{e_k}\tau(\phi)\rangle 
=&0,\\
\langle f\nu,\tau(\phi)\rangle\langle d\phi,\bar\nabla\bar\Delta (|d\phi|^2f\nu)\rangle
=&m^2Hf\langle d\phi,\bar\nabla\bar\Delta (f\nu)\rangle
,\\
\langle f\nu,\tau(\phi)\rangle\langle d\phi,\bar\nabla\big(\langle d\phi,\bar\nabla\tau(\phi)\rangle f\nu\big)\rangle
=&m^4H^4f^2
,\\
\langle f\nu,\tau(\phi)\rangle\langle d\phi,\bar\nabla\big(\langle f\nu,\bar\nabla_{e_k}\tau(\phi)\rangle d\phi(e_k)\big)\rangle=&0,\\
|\bar\nabla_{e_k}\big(\langle f\nu,\tau(\phi)\rangle d\phi(e_k)\big)|^2
=&m^2H^2|\nabla f|^2+m^4H^4f^2,
\\
|\tau(\phi)|^4f^2
=&m^4H^4f^2,\\
|\langle f\nu,\bar\nabla_{e_k}\tau(\phi)\rangle d\phi(e_k)|^2=&0.
\end{align*}
Inserting all the above terms into \eqref{sv-4energy-space-form-hyper} and
also using integration by parts we find
\begin{align*}
Q_4(f\nu,f\nu)=&\int_M\big|\bar\Delta^2(f\nu)-mK\bar\Delta(f\nu)\big|^2\dv \\
&+K\int_M\big(-2m^2H^2|\bar\Delta(f\nu)|^2
+2mH|A|^2\langle\bar\nabla\bar\Delta (f\nu),d\phi\rangle f
+2mH\langle\bar\nabla\bar\Delta^2(f\nu),d\phi\rangle f \\
&-2m^2H^2|A|^2|\nabla f|^2
-8m^2H^2|A|^4f^2
\big)\dv \\
&+K^2\int_M\big(
-2m^2H\langle d\phi,\bar\nabla\bar\Delta(f\nu)\rangle f\\
&+(2m^3H^2+m^2H^2)|\nabla f|^2+(8m^3|A|^2H^2+4m^4H^4)f^2
\big)\dv.
\end{align*}
Some of the above terms need further attention:
Using integration by parts and \eqref{eq:laplace-hypersurface} we find
\begin{align*}
\int_M\langle\bar\nabla\bar\Delta (f\nu),d\phi\rangle f\dv
=-mH\int_M(|\nabla f|^2+|A|^2f^2)\dv
-2\int_M\langle A(\operatorname{grad}f),\operatorname{grad}f\rangle\dv.
\end{align*}
Again, using integration by parts and the self-adjointness of the Laplacian we find
\begin{align*}
\int_M\langle\bar\nabla\bar\Delta^2(f\nu),d\phi\rangle f\dv
=&-mH\int_M|\bar\Delta(f\nu)|^2\dv
-\int_M\langle \bar\Delta^2(f\nu),d\phi(e_j)\rangle\nabla_{e_j}f\dv\\
=&-mH\int_M|\bar\Delta(f\nu)|^2\dv
-2\int_M\langle A(\operatorname{grad}\Delta f),\operatorname{grad}f\rangle\dv \\
&-2|A|^2\int_M\langle A(\operatorname{grad} f),\operatorname{grad} f\rangle\dv
-2\int_M\langle A(\Delta\operatorname{grad} f),\operatorname{grad} f\rangle\dv \\
&-2\int_M\langle A(e_i),A(\operatorname{grad} f)\rangle\langle A(e_i),\operatorname{grad} f\rangle \dv
\end{align*}
completing the proof.
\end{proof}

In order to further manipulate \eqref{qf4-a} we will need the following identities which 
can be obtained by employing \eqref{eq:laplace-hypersurface}, \eqref{eq:nabla-laplace-hypersurface} and
\eqref{eq:laplace-2-hypersurface} such that
\begin{align}
\label{bardelta-identities}
|\bar\Delta(f\nu)|^2=&|\Delta f|^2+|A|^4f^2+2|A|^2f\Delta f+4|A(\operatorname{grad} f)|^2,\\
\nonumber|\bar\nabla\bar\Delta(f\nu)|^2=&|\nabla\Delta f|^2+|A|^4|\nabla f|^2
+4|\langle A(e_i),A(\operatorname{grad}f)\rangle|^2 \\
\nonumber&+2|A|^2\nabla\Delta f\nabla f+4\langle A(\operatorname{grad}\Delta f),A(\operatorname{grad} f)\rangle
+4|A|^2|A(\operatorname{grad}f)|^2 \\
\nonumber&+|A|^2|\Delta f|^2+|A|^6f^2+2|A|^4f\Delta f+4|A(\nabla\operatorname{grad}f)|^2 \\
\nonumber&-4\langle A(e_i),A(\nabla_{e_i}\operatorname{grad}f)\rangle\Delta f
-4|A|^2\langle A(e_i),A(\nabla_{e_i}\operatorname{grad}f)\rangle f, \\
\nonumber|\bar\Delta^2(f\nu)|^2=&|\Delta^2f|^2+4|A|^4|\Delta f|^2+|A|^8f^2
+16|\langle A(e_i),A(\nabla_{e_i}\operatorname{grad} f\rangle|^2 \\
\nonumber&+4|A|^2\Delta^2f\Delta f+2|A|^4f\Delta^2f-8\Delta^2f\langle A(e_i),A(\nabla_{e_i}\operatorname{grad} f\rangle \\
\nonumber&+4|A|^6f\Delta f-16|A|^2\Delta f\langle A(e_i),A(\nabla_{e_i}\operatorname{grad} f)\rangle
-8|A|^4f\langle A(e_i),A(\nabla_{e_i}\operatorname{grad} f)\rangle \\
\nonumber&+4|A(\operatorname{grad}\Delta f)|^2
+4|A|^4|A(\operatorname{grad} f)|^2+4|A(\Delta \operatorname{grad}f)|^2 \\
\nonumber&+4|\langle A(e_i),A(\operatorname{grad}f)\rangle A(e_i)|^2 
+8|A|^2\langle A(\operatorname{grad}\Delta f),A(\operatorname{grad} f)\rangle \\
\nonumber&+8\langle A(\operatorname{grad}\Delta f),A(\Delta \operatorname{grad}f)\rangle
+8\langle A(\operatorname{grad}\Delta f),A(e_i)\rangle\langle A(e_i),A(\operatorname{grad}f)\rangle  \\
\nonumber&+8|A|^2\langle A(\operatorname{grad} f),A(\operatorname{grad} \Delta f)\rangle
+8\langle A(\Delta \operatorname{grad} f),A(e_i)\rangle\langle A(e_i),A(\operatorname{grad}f)\rangle \\
\nonumber&+8|A|^2|\langle A(e_i),A(\nabla_{e_i}\operatorname{grad} f)\rangle|^2.
\end{align}

\subsection{The normal stability of the small 4-harmonic hypersphere}
After having established the necessary formulas for the study of the normal stability
of general 4-harmonic hypersurfaces in Riemannian space forms of constant curvature \(K\)
we now turn to the stability analysis of the small proper 4-harmonic hypersphere
\(\phi\colon\s^m(\frac{1}{2})\to\s^{m+1}\).

The small 4-harmonic hypersphere is characterized by the following geometric data
\begin{align}
\label{data-small-hypersphere}
A=-\sqrt{3}\operatorname{Id},\qquad H=-\sqrt{3},\qquad |A|^2=3m,\qquad H^2=3.
\end{align}

In the following, we will make use of the classic Bochner formulas
\begin{align*}
\int_M|\nabla^2f|^2\dv&=
-\int_M\langle\operatorname{Ric}^M(\nabla f),\nabla f\rangle\dv
+\int_M|\Delta f|^2\dv,\\
\nabla\Delta f&=\Delta\nabla f+\operatorname{Ric}^M(\nabla f),
\end{align*}
which hold on every closed Riemannian manifold \(M\), see for example \cite[Chapter 7]{MR2243772}.

In the case that \(M=\s^m(a)\)
we have \(\operatorname{Ric}^M=\frac{m-1}{a^2}g\) such that
the above formula turns into
\begin{align}
\label{bochner-4hyper}
\int_{\s^m(a)}|\nabla^2f|^2\dv&=-\int_{\s^m(a)}\frac{m-1}{a^2}|\nabla f|^2\dv
+\int_{\s^m(a)}|\Delta f|^2\dv,\\
\nonumber\nabla\Delta f&=\Delta\nabla f+\frac{m-1}{a^2}\nabla f.
\end{align}

With these tools at hand we now give the following 

\begin{Prop}
Let \(\phi\colon \s^m(\frac{1}{2})\to\s^{m+1}\)
be the small proper \(4\)-harmonic hypersphere.
Then, the quadratic form describing its normal stability
is given by
\begin{align}
\label{qf4-b}
Q_4(f\nu,f\nu)
=&\int_{\s^m(\frac{1}{2})}
\big(\big|\bar\Delta^2(f\nu)\big|^2 
-2m|\bar\nabla\bar\Delta(f\nu)|^2
-11m^2|\bar\Delta(f\nu)|^2 \\
\nonumber
&-24m|\Delta f|^2
+(-24m^3-9m^2-84m)|\nabla f|^2-144 m^4f^2 \big)\dv
\end{align}
where \(\nu\) is the unit normal of \(\s^m(\frac{1}{2})\)
and \(f\in C^\infty(\s^m(\frac{1}{2}))\).
\end{Prop}
\begin{proof}
Setting \(K=1\) in \eqref{qf4-a} we find
\begin{align*}
Q_4(f\nu,f\nu)=
&\int_{\s^m(\frac{1}{2})}\big|\bar\Delta^2(f\nu)-m\bar\Delta(f\nu)\big|^2\dv \\
\nonumber&+\int_{\s^m(\frac{1}{2})}
\bigg(
-4m^2H^2|\bar\Delta(f\nu)|^2 
-4m^2H^2|A|^2|\nabla f|^2 
-10m^2H^2|A|^4f^2
\\
\nonumber&
-4mH\langle A(\operatorname{grad}\Delta f),\operatorname{grad}f\rangle 
-4mH\langle A(e_i),A(\operatorname{grad} f)\rangle\langle A(e_i),\operatorname{grad} f\rangle
\\
\nonumber&-8mH|A|^2\langle A(\operatorname{grad} f),\operatorname{grad} f\rangle 
-4mH\langle A(\Delta\operatorname{grad} f),\operatorname{grad} f\rangle 
\nonumber 
\\
&+(4m^3H^2+m^2H^2)|\nabla f|^2+(10m^3|A|^2H^2+4m^4H^4)f^2 \\
\nonumber&+ 4m^2H\langle A(\operatorname{grad} f),\operatorname{grad} f\rangle
\bigg)\dv \\
=&\int_{\s^m(\frac{1}{2})}
\big(\big|\bar\Delta^2(f\nu)\big|^2
-2m|\bar\nabla\bar\Delta(f\nu)|^2
-11m^2|\bar\Delta(f\nu)|^2 \\
&-12m|\Delta f|^2-12m\Delta\nabla f\nabla f \\
&+(-24m^3-57m^2-36m)|\nabla f|^2-144 m^4f^2 \big)\dv,
\end{align*}
where we have used \eqref{data-small-hypersphere} in the second step. 
The claim now follows from interchanging derivatives
and using integration by parts.
\end{proof}

We are now ready to prove the first main result of this article:
\begin{Satz}
\label{thm:4hyper-stab}
Let \(\phi\colon\s^m(\frac{1}{2})\to\s^{m+1}\)
be the small proper 4-harmonic hypersphere. Then, the quadratic form 
describing its normal stability is given by
\begin{align}
\label{qf4-c}
Q_4(f\nu,f\nu)=\int_{\s^m(\frac{1}{2})}
\big(&|\Delta^2f|^2+(10m+72)|\nabla\Delta f|^2\\
\nonumber&+(25m^2+286m+480)|\Delta f|^2 \\
\nonumber&+(-36m^3+327m^2-924m+588)|\nabla f|^2 \\
\nonumber&-216m^4f^2
\big)\dv.
\end{align}
In particular, the normal index of the small proper 4-harmonic hypersphere
is equal to one, i.e.
\begin{align*}
\operatorname{Ind}^{\rm{nor}}_4(\s^m(\frac{1}{2})\to\s^{m+1})=1.
\end{align*}
\end{Satz}

\begin{proof}
First, using the identities provided by \eqref{bardelta-identities} and
inserting \eqref{data-small-hypersphere} we find
\begin{align*}
\int_{\s^m(\frac{1}{2})}|\bar\Delta(f\nu)|^2\dv=
\int_{\s^m(\frac{1}{2})}\big(|\Delta f|^2+(6m+12)|\nabla f|^2+9m^2f^2\big)\dv
\end{align*}
and also
\begin{align*}
\int_{\s^m(\frac{1}{2})}|\bar\nabla\bar\Delta(f\nu)|^2\dv=
\int_{\s^m(\frac{1}{2})}\big(&|\nabla\Delta f|^2+(9m+24)|\Delta f|^2 
+(27m^2+72m+36)|\nabla f|^2 \\
&+27m^3f^2
+12|\nabla\nabla f|^2
\big)\dv\\
=\int_{\s^m(\frac{1}{2})}\big(&|\nabla\Delta f|^2
+(9m+36)|\Delta f|^2+(27m^2+24m+84)|\nabla f|^2 \\
&+27m^3f^2\big)\dv.
\end{align*}
Note that we used the Bochner formulas \eqref{bochner-4hyper} to interchange covariant derivatives 
several times.

Inserting the data \eqref{data-small-hypersphere} into the third equation of
\eqref{bardelta-identities} we find
\begin{align*}
\nonumber|\bar\Delta^2(f\nu)|^2=&|\Delta^2f|^2+36m^2|\Delta f|^2+81m^4f^2
+144|\Delta f|^2 \\
\nonumber&+12m\Delta^2f\Delta f+18m^2f\Delta^2f+24\Delta^2f\Delta f \\
\nonumber&+108m^3f\Delta f+144m|\Delta f|^2
+216m^2f\Delta f \\
\nonumber&+12|\nabla\Delta f|^2
+108m^2|\nabla f|^2+12|\Delta\nabla f|^2 \\
\nonumber&+108|\nabla f|^2 
+72m\nabla\Delta f\nabla f \\
\nonumber&+24\nabla\Delta f\Delta\nabla f
+72\nabla\Delta f\nabla f  \\
\nonumber&+72m\nabla f\nabla\Delta f
+72\Delta\nabla f\nabla f \\
\nonumber&+216m|\Delta f|^2.    
\end{align*}

This allows us to conclude that
\begin{align*}
\int_{\s^m(\frac{1}{2})}|\bar\Delta^2(f\nu)|^2\dv=
\int_{\s^m(\frac{1}{2})}
\big(&
|\Delta^2f|^2 +(12m+36)|\nabla\Delta f|^2
+(54m^2+504m+216)|\Delta f|^2 \\
&+(108m^3+324m^2+108)|\nabla f|^2
+81m^4f^2 \\
&+12|\Delta\nabla f|^2+24\nabla\Delta f\Delta\nabla f+72\Delta\nabla f\nabla f
\big)\dv\\
=\int_{\s^m(\frac{1}{2})}
\big(&
|\Delta^2f|^2+(12m+72)|\nabla\Delta f|^2+(54m^2+382m+480)|\Delta f|^2 \\
&+(108m^3+516m^2-672m+588)|\nabla f|^2+81m^4f^2
\big)\dv,    
\end{align*}
where we again made use of the Bochner formulas \eqref{bochner-4hyper}.

Combining the above identities with \eqref{qf4-b} now completes the proof of the first part of the theorem.

Regarding the second claim concerning the normal index of the small proper 4-harmonic hypersphere
we note that the quadratic form \eqref{qf4-c} will be negative for \(f=const\)
which corresponds to an eigenfunction of the Laplace operator with eigenvalue \(0\).
Now, let \(\lambda\) be a non-zero eigenvalue of the Laplace operator, 
that is \(\Delta f=\lambda f\).
Then, the quadratic form \eqref{qf4-c} simplifies to
\begin{align}
\label{qf4-d}
Q_4(f\nu,f\nu)=\int_{\s^m(\frac{1}{2})}
\big(&\lambda^4+(10m+72)\lambda^3\\
\nonumber&+(25m^2+286m+480)\lambda^2 \\
\nonumber&+(-36m^3+327m^2-924m+588)\lambda \\
\nonumber&-216m^4
\big)f^2\dv.
\end{align}
The first non-zero eigenvalue of the Laplace operator on \(\s^m(\frac{1}{2})\)
is \(\lambda_1=4m\) in which case the quadratic form \eqref{qf4-d} acquires the form
\begin{align*}
Q_4(f\nu,f\nu)=\int_{\s^m(\frac{1}{2})}
\big(
936m^4+10672m^3+3984m^2+2352m
\big)f^2\dv.
\end{align*}
As this expression is clearly positive for all \(m\) the second claim of the theorem follows directly.
Also, it is not difficult to see that for all higher eigenvalues of the Laplacian the quadratic form \eqref{qf4-d}
will be positive.
\end{proof}

\section{The normal stability of the ES-4 curvature term}
In this section we investigate the influence of the curvature term appearing in the ES-4-energy 
\eqref{energy-es4-harmonic}
on the normal stability of its critical points.
First, recall that the curvature term of the ES-4-energy is given by
\begin{align*}
\hat E_4(\phi)=\frac{1}{4}\int_M|R^N(d\phi(e_i),d\phi(e_j))\tau(\phi)|^2\dv 
\end{align*}
and its first variation formula (for the full ES-4-tension field see \eqref{es-4-tension}) acquires the form

\begin{align}
\label{first-variation-ct}
\frac{d}{dt}\frac{1}{4}\int_M|R^N(d\phi_t(e_i),d\phi_t(e_j))\tau(\phi_t)|^2\dv 
=-\int_M\langle d\phi_t(\partial_t),\hat\tau_4(\phi_t)\rangle\dv,
\end{align}
where \(\phi_t\) represents a variation of the map \(\phi\) as defined in \eqref{dfn:variation-phi}.

Here, \(\hat\tau_4(\phi)\) is the quantity defined by
\begin{align*}
\hat{\tau}_4(\phi)=-\frac{1}{2}\big(2\xi_1+2d^\ast\Omega_1+\bar\Delta\Omega_0+\tr R^N(d\phi(\cdot),\Omega_0)d\phi(\cdot)\big),
\end{align*}
where we have used the following abbreviations
\begin{align*}
\Omega_0&=R^N(d\phi(e_i),d\phi(e_j))R^N(d\phi(e_i),d\phi(e_j))\tau(\phi), \\
\nonumber\Omega_1(X)&=R^N(R^N(d\phi(X),d\phi(e_j))\tau(\phi),\tau(\phi))d\phi(e_j),\\
\nonumber\xi_1&=-(\nabla_{d\phi(e_j)}R^N)(R^N(d\phi(e_i),d\phi(e_j))\tau(\phi),\tau(\phi))d\phi(e_i).
\end{align*}

For a detailed derivation of the first variation formula \eqref{first-variation-ct} we refer to \cite[Section 3]{MR4106647}.

\begin{Lem}[Second variation of the ES-4-curvature term]
Let \(\phi\colon M\to N\) be a smooth map and consider a variation of the map \(\phi\)
as described in \eqref{dfn:variation-phi}. Then, the following formula holds:
\begin{align}
\label{sv-es4-curvature}
\frac{d^2}{dt^2}&\frac{1}{4}\int_M|R^N(d\phi_t(e_i),d\phi_t(e_j))\tau(\phi_t)|^2\dv\big|_{t=0} \\
\nonumber=&-\int_M\langle \frac{\bar\nabla}{\partial t}d\phi_t(\partial_t),\hat{\tau}_4(\phi_t)
\rangle\big|_{t=0}\dv \\
\nonumber&
-\int_M \langle V,(\nabla_{V}\nabla_{d\phi(e_j)}R^N)(R^N(d\phi(e_i),d\phi(e_j))\tau(\phi),\tau(\phi))d\phi(e_i)\rangle\dv \\
\nonumber&-\int_M\langle V,(\nabla_{\bar\nabla_{e_j}V}R^N)(R^N(d\phi(e_i),d\phi(e_j))\tau(\phi),\tau(\phi))d\phi(e_i)\rangle\dv
\\
\nonumber&-\int_M\langle V,\nabla_{d\phi(e_j)}R^N)((\nabla_{V}R^N)(d\phi(e_i),d\phi(e_j))\tau(\phi),\tau(\phi))d\phi(e_i)
\rangle\dv\\
\nonumber&-\int_M\langle V,(\nabla_{d\phi(e_j)}R^N)(R^N(\bar\nabla_{e_i}V,d\phi(e_j))\tau(\phi),\tau(\phi))d\phi(e_i)\rangle\dv \\
\nonumber&-\int_M\langle V,(\nabla_{d\phi(e_j)}R^N)(R^N(d\phi(e_i),\bar\nabla_{e_j}V)\tau(\phi),\tau(\phi))d\phi(e_i)\rangle\dv \\
\nonumber&+\int_M\langle V,(\nabla_{d\phi(e_j)}R^N)(R^N(d\phi(e_i),d\phi(e_j))\bar\Delta V,\tau(\phi))d\phi(e_i)\rangle\dv \\
\nonumber&-\int_M\langle V,(\nabla_{d\phi(e_j)}R^N)(R^N(d\phi(e_i),d\phi(e_j))R^N(V,d\phi(e_k))d\phi(e_k),\tau(\phi))d\phi(e_i)\rangle\dv \\
\nonumber&+\int_M\langle V,(\nabla_{d\phi(e_j)}R^N)(R^N(d\phi(e_i),d\phi(e_j))\tau(\phi),\bar\Delta V)d\phi(e_i)\rangle\dv \\
\nonumber&-\int_M\langle V,(\nabla_{d\phi(e_j)}R^N)(R^N(d\phi(e_i),d\phi(e_j))\tau(\phi),R^N(V,d\phi(e_k))d\phi(e_k))d\phi(e_i) \rangle\dv\\
\nonumber&-\int_M\langle V,(\nabla_{d\phi(e_j)}R^N)(R^N(d\phi(e_i),d\phi(e_j))\tau(\phi),\tau(\phi))\bar\nabla_{e_i}V\rangle\dv
\\
\nonumber&-\int_M\langle V,R^N(V,d\phi(e_i))R^N(R^N(d\phi(e_i),d\phi(e_j))\tau(\phi),\tau(\phi))d\phi(e_j)
\rangle\dv \\
\nonumber&+\int_M\langle \bar\nabla_{e_i}V,(\nabla_{V}R^N)(R^N(d\phi(e_i),d\phi(e_j))\tau(\phi),\tau(\phi))d\phi(e_j)\rangle\dv
\\
\nonumber&+\int_M\langle \bar\nabla_{e_i}V,R^N((\nabla_{V}R^N)(d\phi(e_i),d\phi(e_j))\tau(\phi),\tau(\phi))d\phi(e_j)\rangle\dv
\\
\nonumber&+\int_M\langle \bar\nabla_{e_i}V,R^N(R^N(\bar\nabla_{e_i}V,d\phi(e_j))\tau(\phi),\tau(\phi))d\phi(e_j)\rangle\dv
\\
\nonumber&+\int_M\langle \bar\nabla_{e_i}V,R^N(R^N(d\phi(e_i),\bar\nabla_{e_j}V)\tau(\phi),\tau(\phi))d\phi(e_j)\rangle\dv
\\
\nonumber&-\int_M\langle \bar\nabla_{e_i}V,R^N(R^N(d\phi(e_i),d\phi(e_j))\bar\Delta V,\tau(\phi))d\phi(e_j)
\rangle\dv
\\
\nonumber&+\int_M\langle \bar\nabla_{e_i}V,R^N(R^N(d\phi(e_i),d\phi(e_j))
R^N(V,d\phi(e_k))d\phi(e_k),\tau(\phi))d\phi(e_j)\rangle\dv
\\
\nonumber&-\int_M\langle \bar\nabla_{e_i}V,R^N(R^N(d\phi(e_i),d\phi(e_j))\tau(\phi),\bar\Delta V)d\phi(e_j)\rangle\dv
\\
\nonumber&+\int_M\langle \bar\nabla_{e_i}V,R^N(R^N(d\phi(e_i),d\phi(e_j))\tau(\phi),R^N(V,d\phi(e_k))d\phi(e_k))
d\phi(e_j)\rangle\dv
\\
\nonumber&+\int_M\langle \bar\nabla_{e_i}V,R^N(R^N(d\phi(e_i),d\phi(e_j))\tau(\phi),\tau(\phi))\bar\nabla_{e_j}V\rangle\dv\\
\nonumber&
-\frac{1}{2}\int_M\langle V,R^N(V,d\phi(e_k))\bar\nabla_{e_k}\big(
R^N(d\phi(e_i),d\phi(e_j))R^N(d\phi(e_i),d\phi(e_j))\tau(\phi)
\big)\rangle\dv\\
\nonumber&
+\frac{1}{2}\int_M\langle\bar\nabla_{e_k}V,
R^N(V,d\phi(e_k))R^N(d\phi(e_i),d\phi(e_j))R^N(d\phi(e_i),d\phi(e_j))\tau(\phi)\rangle\dv
\\
\nonumber&+\frac{1}{2}\int_M\langle\bar\Delta V,(\nabla_{V}R^N)(d\phi(e_i),d\phi(e_j))R^N(d\phi(e_i),d\phi(e_j))\tau(\phi)\rangle\dv \\
\nonumber&+\int_M\langle\bar\Delta V,R^N(\bar\nabla_{e_i}V,d\phi(e_j))R^N(d\phi(e_i),d\phi(e_j))\tau(\phi)\rangle\dv \\
\nonumber&+\frac{1}{2}\int_M\langle\bar\Delta V,R^N(d\phi(e_i),d\phi(e_j))(\nabla_{V}R^N)(d\phi(e_i),d\phi(e_j))\tau(\phi)\rangle\dv \\
\nonumber&+\int_M\langle\bar\Delta V,R^N(d\phi(e_i),d\phi(e_j))R^N(\bar\nabla_{e_i}V,d\phi(e_j))\tau(\phi)\rangle\dv \\
\nonumber&-\frac{1}{2}\int_M\langle\bar\Delta V,R^N(d\phi(e_i),d\phi(e_j))R^N(d\phi(e_i),d\phi(e_j))\bar\Delta V\rangle\dv \\
\nonumber&+\frac{1}{2}\int_M\langle\bar\Delta V,R^N(d\phi(e_i),d\phi(e_j))R^N(d\phi(e_i),d\phi(e_j))R^N(V,d\phi(e_k))d\phi(e_k)\rangle\dv\\
\nonumber&+\frac{1}{2}\int_M\langle V,(\nabla_{V}R^N)\big(d\phi(e_k),R^N(d\phi(e_i),d\phi(e_j))R^N(d\phi(e_i),d\phi  (e_j))\tau(\phi)\big)d\phi(e_k)\rangle\dv\\
\nonumber&+\frac{1}{2}\int_M\langle V,R^N\big(\bar\nabla_{e_k}V,R^N(d\phi(e_i),d\phi(e_j))R^N(d\phi(e_i),d\phi(e_j))\tau(\phi)\big)d\phi(e_k)\rangle\dv\\
\nonumber&+\frac{1}{2}\int_M\langle V,R^N\big(d\phi(e_k),(\nabla_{V}R^N)(d\phi(e_i),d\phi(e_j))R^N(d\phi(e_i),d\phi  (e_j))\tau(\phi)\big)d\phi(e_k)\rangle\dv\\
\nonumber&+\int_M\langle V,R^N\big(d\phi(e_k),R^N(\bar\nabla_{e_i}V,d\phi(e_j))R^N(d\phi(e_i),d\phi(e_j))\tau(\phi)\big)d\phi(e_k)\rangle\dv\\
\nonumber&+\frac{1}{2}\int_M\langle V,R^N\big(d\phi(e_k),R^N(d\phi(e_i),d\phi(e_j))(\nabla_{V}R^N)(d\phi(e_i),d\phi (e_j))\tau(\phi)\big)d\phi(e_k)\rangle\dv\\
\nonumber&+\int_M\langle V,R^N\big(d\phi(e_k),R^N(d\phi(e_i),d\phi(e_j))R^N(\bar\nabla_{e_i}V,d\phi(e_j))\tau(\phi)\big)d\phi(e_k)\rangle\dv\\
\nonumber&-\frac{1}{2}\int_M\langle V,R^N\big(d\phi(e_k),R^N(d\phi(e_i),d\phi(e_j))R^N(d\phi(e_i),d\phi(e_j))\bar\Delta 
V\big)d\phi(e_k)\rangle\dv\\
\nonumber&+\frac{1}{2}\int_M\langle V,R^N\big(d\phi(e_k),R^N(d\phi(e_i),d\phi(e_j))R^N(d\phi(e_i),d\phi(e_j))
R^N(V,d\phi(e_l))d\phi(e_l)\big)d\phi(e_k)\rangle\dv\\
\nonumber&+\frac{1}{2}\int_M\langle V,R^N\big(d\phi(e_k),R^N(d\phi(e_i),d\phi(e_j))R^N(d\phi(e_i),d\phi(e_j))\tau(\phi)\big)\bar\nabla_{e_k}V\rangle\dv.
\end{align}

\end{Lem}
\begin{proof}
We differentiate \eqref{first-variation-ct} once more and find
\begin{align*}
\frac{d^2}{dt^2}\frac{1}{4}&\int_M|R^N(d\phi_t(e_i),d\phi_t(e_j))\tau(\phi_t)|^2\dv \\
\nonumber=&-\int_M\langle \frac{\bar\nabla}{\partial t}d\phi_t(\partial_t),\hat\tau_4(\phi_t)\rangle\dv \\
\nonumber&+\int_M\langle d\phi_t(\partial_t),\frac{\bar\nabla}{\partial t}\xi_1\rangle\dv
+\int_M\langle d\phi_t(\partial_t),\frac{\bar\nabla}{\partial t}d^\ast\Omega_1\rangle\dv
+\frac{1}{2}\int_M\langle d\phi_t(\partial_t),\frac{\bar\nabla}{\partial t}\bar\Delta\Omega_0\rangle \dv \\
\nonumber&+\frac{1}{2}
\int_M\langle d\phi_t(\partial_t),\frac{\bar\nabla}{\partial t}
\big(\tr R^N(d\phi_t(\cdot),\Omega_0)d\phi_t(\cdot)\big)\rangle\dv.
\end{align*}

First, we manipulate the term containing \(\xi_1\).

\begin{align*}
\frac{\bar\nabla}{\partial t}\xi_1
=&\frac{\bar\nabla}{\partial t}\big(-(\nabla_{d\phi_t(e_j)}R^N)(R^N(d\phi_t(e_i),d\phi_t(e_j))\tau(\phi_t),\tau(\phi_t))d\phi_t(e_i)\big) \\
=&-(\nabla_{d\phi_t(\partial_t)}\nabla_{d\phi_t(e_j)}R^N)(R^N(d\phi_t(e_i),d\phi_t(e_j))\tau(\phi_t),\tau(\phi_t))d\phi_t(e_i) \\
&-(\nabla_{\bar\nabla_{e_j}d\phi_t(\partial_t)}R^N)(R^N(d\phi_t(e_i),d\phi_t(e_j))\tau(\phi_t),\tau(\phi_t))d\phi_t(e_i)
\\
&-(\nabla_{d\phi_t(e_j)}R^N)((\nabla_{d\phi_t(\partial_t)}R^N)(d\phi_t(e_i),d\phi_t(e_j))\tau(\phi_t),\tau(\phi_t))d\phi_t(e_i) \\
&-(\nabla_{d\phi_t(e_j)}R^N)(R^N(\bar\nabla_{e_i}d\phi_t(\partial_t),d\phi_t(e_j))\tau(\phi_t),\tau(\phi_t))d\phi_t(e_i) \\
&-(\nabla_{d\phi_t(e_j)}R^N)(R^N(d\phi_t(e_i),\bar\nabla_{e_j}d\phi_t(\partial_t))\tau(\phi_t),\tau(\phi_t))d\phi_t(e_i) \\
&+(\nabla_{d\phi_t(e_j)}R^N)(R^N(d\phi_t(e_i),d\phi_t(e_j))\bar\Delta d\phi_t(\partial_t),\tau(\phi_t))d\phi_t(e_i) \\
&-(\nabla_{d\phi_t(e_j)}R^N)(R^N(d\phi_t(e_i),d\phi_t(e_j))R^N(d\phi_t(\partial_t),d\phi_t(e_k))d\phi_t(e_k),\tau(\phi_t))d\phi_t(e_i) \\
&+(\nabla_{d\phi_t(e_j)}R^N)(R^N(d\phi_t(e_i),d\phi_t(e_j))\tau(\phi_t),\bar\Delta d\phi_t(\partial_t))d\phi_t(e_i) \\
&-(\nabla_{d\phi_t(e_j)}R^N)(R^N(d\phi_t(e_i),d\phi_t(e_j))\tau(\phi_t),R^N(d\phi_t(\partial_t),d\phi_t(e_k))d\phi_t(e_k))d\phi_t(e_i) \\
&-(\nabla_{d\phi_t(e_j)}R^N)(R^N(d\phi_t(e_i),d\phi_t(e_j))\tau(\phi_t),\tau(\phi_t))\bar\nabla_{e_i}d\phi_t(\partial_t).
\end{align*}
Evaluating at \(t=0\) then gives
\begin{align*}
&\int_M\langle d\phi_t(\partial_t),\frac{\bar\nabla}{\partial t}\xi_1\rangle\dv \big|_{t=0}\\
=&-\int_M \langle V,(\nabla_{V}\nabla_{d\phi(e_j)}R^N)(R^N(d\phi(e_i),d\phi(e_j))\tau(\phi),\tau(\phi))d\phi(e_i)\rangle\dv \\
&-\int_M\langle V,(\nabla_{\bar\nabla_{e_j}V}R^N)(R^N(d\phi(e_i),d\phi(e_j))\tau(\phi),\tau(\phi))d\phi(e_i)\rangle\dv
\\
&-\int_M\langle V,\nabla_{d\phi(e_j)}R^N)((\nabla_{V}R^N)(d\phi(e_i),d\phi(e_j))\tau(\phi),\tau(\phi))d\phi(e_i)
\rangle\dv\\
&-\int_M\langle V,(\nabla_{d\phi(e_j)}R^N)(R^N(\bar\nabla_{e_i}V,d\phi(e_j))\tau(\phi),\tau(\phi))d\phi(e_i)\rangle\dv \\
&-\int_M\langle V,(\nabla_{d\phi(e_j)}R^N)(R^N(d\phi(e_i),\bar\nabla_{e_j}V)\tau(\phi),\tau(\phi))d\phi(e_i)\rangle\dv \\
&+\int_M\langle V,(\nabla_{d\phi(e_j)}R^N)(R^N(d\phi(e_i),d\phi(e_j))\bar\Delta V,\tau(\phi))d\phi(e_i)\rangle\dv \\
&-\int_M\langle V,(\nabla_{d\phi(e_j)}R^N)(R^N(d\phi(e_i),d\phi(e_j))R^N(V,d\phi(e_k))d\phi(e_k),\tau(\phi))d\phi(e_i)\rangle\dv \\
&+\int_M\langle V,(\nabla_{d\phi(e_j)}R^N)(R^N(d\phi(e_i),d\phi(e_j))\tau(\phi),\bar\Delta V)d\phi(e_i)\rangle\dv \\
&-\int_M\langle V,(\nabla_{d\phi(e_j)}R^N)(R^N(d\phi(e_i),d\phi(e_j))\tau(\phi),R^N(V,d\phi(e_k))d\phi(e_k))d\phi(e_i) \rangle\dv\\
&-\int_M\langle V,(\nabla_{d\phi(e_j)}R^N)(R^N(d\phi(e_i),d\phi(e_j))\tau(\phi),\tau(\phi))\bar\nabla_{e_i}V\rangle\dv.
\end{align*}

In order to manipulate the term containing the divergence of \(\Omega_1\) we calculate
\begin{align*}
\int_M\langle d\phi_t(\partial_t),\frac{\bar\nabla}{\partial t}d^\ast\Omega_1\rangle\big|_{t=0}\dv
=&-\int_M\langle d\phi_t(\partial_t),\frac{\bar\nabla}{\partial t}\bar\nabla_{e_i}\Omega_1(e_i)
\rangle\big|_{t=0}\dv \\
=&-\int_M\langle d\phi_t(\partial_t),R^N(d\phi_t(\partial_t),d\phi_t(e_i))\Omega_1(e_i)
\rangle\big|_{t=0}\dv \\
&+\int_M\langle \bar\nabla_{e_i}d\phi_t(\partial_t),\frac{\bar\nabla}{\partial t}\Omega_1(e_i)
\rangle\big|_{t=0}\dv,
\end{align*}
where we first interchanged covariant derivatives and used integration by parts in the second step.
Now, a direct calculation using \eqref{commutator-t-tension} shows that
\begin{align*}
\frac{\bar\nabla}{\partial t}\Omega_1(e_i)=&
\frac{\bar\nabla}{\partial t}\big(R^N(R^N(d\phi_t(e_i),d\phi_t(e_j))\tau(\phi_t),\tau(\phi_t))d\phi_t(e_j)\big)\\
=&(\nabla_{d\phi_t(\partial_t)}R^N)(R^N(d\phi_t(e_i),d\phi_t(e_j))\tau(\phi_t),\tau(\phi_t))d\phi_t(e_j)\\
&+R^N((\nabla_{d\phi_t(\partial_t)}R^N)(d\phi_t(e_i),d\phi_t(e_j))\tau(\phi_t),\tau(\phi_t))d\phi_t(e_j)\\
&+R^N(R^N(\bar\nabla_{e_i}d\phi_t(\partial_t),d\phi_t(e_j))\tau(\phi_t),\tau(\phi_t))d\phi_t(e_j) \\
&+R^N(R^N(d\phi_t(e_i),\bar\nabla_{e_j}d\phi_t(\partial_t))\tau(\phi_t),\tau(\phi_t))d\phi(e_j) \\
&-R^N(R^N(d\phi_t(e_i),d\phi_t(e_j))\bar\Delta d\phi_t(\partial_t),\tau(\phi_t))d\phi(e_j) \\
&+R^N(R^N(d\phi_t(e_i),d\phi_t(e_j))
R^N(d\phi_t(\partial_t),d\phi_t(e_k))d\phi_t(e_k),\tau(\phi_t))d\phi_t(e_j) \\
&-R^N(R^N(d\phi_t(e_i),d\phi_t(e_j))\tau(\phi_t),\bar\Delta d\phi_t(\partial_t))d\phi_t(e_j) \\
&+R^N(R^N(d\phi_t(e_i),d\phi_t(e_j))\tau(\phi_t),R^N(d\phi_t(\partial_t),d\phi_t(e_k))d\phi_t(e_k))
d\phi_t(e_j) \\
&+R^N(R^N(d\phi_t(e_i),d\phi_t(e_j))\tau(\phi_t),\tau(\phi_t))\bar\nabla_{e_j}d\phi_t(\partial_t). 
\end{align*}
Hence, we can deduce that
\begin{align*}
\int_M&\langle d\phi_t(\partial_t),\frac{\bar\nabla}{\partial t}d^\ast\Omega_1\rangle\big|_{t=0}\dv 
\\
=&-\int_M\langle V,R^N(V,d\phi(e_i))\Omega_1(e_i)
\rangle\dv \\
&+\int_M\langle \bar\nabla_{e_i}V,(\nabla_{V}R^N)(R^N(d\phi(e_i),d\phi(e_j))\tau(\phi),\tau(\phi))d\phi(e_j)\rangle\dv
\\
&+\int_M\langle \bar\nabla_{e_i}V,R^N((\nabla_{V}R^N)(d\phi(e_i),d\phi(e_j))\tau(\phi),\tau(\phi))d\phi(e_j)\rangle\dv
\\
&+\int_M\langle \bar\nabla_{e_i}V,R^N(R^N(\bar\nabla_{e_i}V,d\phi(e_j))\tau(\phi),\tau(\phi))d\phi(e_j)\rangle\dv
\\
&+\int_M\langle \bar\nabla_{e_i}V,R^N(R^N(d\phi(e_i),\bar\nabla_{e_j}V)\tau(\phi),\tau(\phi))d\phi(e_j)\rangle\dv
\\
&-\int_M\langle \bar\nabla_{e_i}V,R^N(R^N(d\phi(e_i),d\phi(e_j))\bar\Delta V,\tau(\phi))d\phi(e_j)
\rangle\dv
\\
&+\int_M\langle \bar\nabla_{e_i}V,R^N(R^N(d\phi(e_i),d\phi(e_j))
R^N(V,d\phi(e_k))d\phi(e_k),\tau(\phi))d\phi(e_j)\rangle\dv
\\
&-\int_M\langle \bar\nabla_{e_i}V,R^N(R^N(d\phi(e_i),d\phi(e_j))\tau(\phi),\bar\Delta V)d\phi(e_j)\rangle\dv
\\
&+\int_M\langle \bar\nabla_{e_i}V,R^N(R^N(d\phi(e_i),d\phi(e_j))\tau(\phi),R^N(V,d\phi(e_k))d\phi(e_k))
d\phi(e_j)\rangle\dv
\\
&+\int_M\langle \bar\nabla_{e_i}V,R^N(R^N(d\phi(e_i),d\phi(e_j))\tau(\phi),\tau(\phi))\bar\nabla_{e_j}V\rangle\dv.
\end{align*}

Concerning the term involving \(\bar\Delta\Omega_0\) we find
\begin{align*}
\int_M\langle d\phi_t(\partial_t),\frac{\bar\nabla}{\partial t}\bar\Delta\Omega_0\rangle\dv
\big|_{t=0}
=\int_M\langle d\phi_t(\partial_t),[\frac{\bar\nabla}{\partial t},\bar\Delta]\Omega_0\rangle\dv
\big|_{t=0}
+\int_M\langle\bar\Delta d\phi_t(\partial_t),\frac{\bar\nabla}{\partial t}\Omega_0\rangle\dv
\big|_{t=0}.
\end{align*}

For any \(X\in\Gamma(\phi^\ast TN)\) it is straightforward to derive that
\begin{align*}
[\frac{\bar\nabla}{\partial t},\bar\Delta]X=
-R^N(d\phi_t(\partial_t),d\phi_t(e_j))\bar\nabla_{e_j}X
-\bar\nabla_{e_j}\big(R^N(d\phi_t(\partial_t),d\phi_t(e_j))X\big)
\end{align*}
such that we find

\begin{align*}
\int_M\langle d\phi_t(\partial_t),[\frac{\bar\nabla}{\partial t},\bar\Delta]\Omega_0\rangle\dv
\big|_{t=0}
=&-\int_M\langle V,R^N(V,d\phi(e_k))\bar\nabla_{e_k}\Omega_0\rangle\dv \\
&-\int_M\langle V,\bar\nabla_{e_k}\big(R^N(V,d\phi(e_k))\Omega_0\big)\rangle\dv.
\end{align*}

Furthermore, a direct calculation shows that
\begin{align*}
\frac{\bar\nabla}{\partial t}\Omega_0=&
\frac{\bar\nabla}{\partial t}\big(R^N(d\phi_t(e_i),d\phi_t(e_j))R^N(d\phi_t(e_i),d\phi_t(e_j))\tau(\phi_t)\big)\\
=&(\nabla_{d\phi_t(\partial_t)}R^N)(d\phi_t(e_i),d\phi_t(e_j))R^N(d\phi_t(e_i),d\phi_t(e_j))\tau(\phi_t) \\
&+2R^N(\bar\nabla_{e_i}d\phi_t(\partial_t),d\phi_t(e_j))R^N(d\phi_t(e_i),d\phi_t(e_j))\tau(\phi_t)\\
&+R^N(d\phi_t(e_i),d\phi_t(e_j))(\nabla_{d\phi_t(\partial_t)}R^N)(d\phi_t(e_i),d\phi_t(e_j))\tau(\phi_t) \\
&+2R^N(d\phi_t(e_i),d\phi_t(e_j))R^N(\bar\nabla_{e_i}d\phi_t(\partial_t),d\phi_t(e_j))\tau(\phi_t)\\
&-R^N(d\phi_t(e_i),d\phi_t(e_j))R^N(d\phi_t(e_i),d\phi_t(e_j))\bar\Delta d\phi_t(\partial_t)\\
&+R^N(d\phi_t(e_i),d\phi_t(e_j))R^N(d\phi_t(e_i),d\phi_t(e_j))R^N(d\phi_t(\partial_t),d\phi_t(e_k))d\phi_t(e_k),
\end{align*}
where we used \eqref{commutator-t-tension} for the last two terms.
Hence, we may conclude that
\begin{align*}
\int_M\langle&\bar\Delta d\phi_t(\partial_t),\frac{\bar\nabla}{\partial t}\Omega_0\rangle\dv
\big|_{t=0}\\
=&\int_M\langle\bar\Delta V,(\nabla_{V}R^N)(d\phi(e_i),d\phi(e_j))R^N(d\phi(e_i),d\phi(e_j))\tau(\phi)\rangle\dv \\
&+2\int_M\langle\bar\Delta V,R^N(\bar\nabla_{e_i}V,d\phi(e_j))R^N(d\phi(e_i),d\phi(e_j))\tau(\phi)\rangle\dv \\
&+\int_M\langle\bar\Delta V,R^N(d\phi(e_i),d\phi(e_j))(\nabla_{V}R^N)(d\phi(e_i),d\phi(e_j))\tau(\phi)\rangle\dv \\
&+2\int_M\langle\bar\Delta V,R^N(d\phi(e_i),d\phi(e_j))R^N(\bar\nabla_{e_i}V,d\phi(e_j))\tau(\phi)\rangle\dv \\
&-\int_M\langle\bar\Delta V,R^N(d\phi(e_i),d\phi(e_j))R^N(d\phi(e_i),d\phi(e_j))\bar\Delta V\rangle\dv \\
&+\int_M\langle\bar\Delta V,R^N(d\phi(e_i),d\phi(e_j))R^N(d\phi(e_i),d\phi(e_j))R^N(V,d\phi(e_k))d\phi(e_k)\rangle.
\end{align*}

Finally, we calculate
\begin{align*}
\frac{\bar\nabla}{\partial t}&\big(\tr R^N(d\phi_t(\cdot),\Omega_0)d\phi_t(\cdot)\big) \\
=&\frac{\bar\nabla}{\partial t}\bigg(R^N\big(d\phi_t(e_k),R^N(d\phi_t(e_i),d\phi_t(e_j))R^N(d\phi_t(e_i),d\phi_t(e_j))\tau(\phi_t)\big)d\phi_t(e_k)\bigg) \\
=&
(\nabla_{d\phi_t(\partial_t)}R^N)\big(d\phi_t(e_k),R^N(d\phi_t(e_i),d\phi_t(e_j))R^N(d\phi_t(e_i),d\phi_t(e_j))\tau(\phi_t)\big)d\phi_t(e_k)\\
&+R^N\big(\bar\nabla_{e_k}d\phi_t(\partial_t),R^N(d\phi_t(e_i),d\phi_t(e_j))R^N(d\phi_t(e_i),d\phi_t(e_j))\tau(\phi_t)\big)d\phi_t(e_k)\\
&+R^N\big(d\phi_t(e_k),(\nabla_{d\phi_t(\partial_t)}R^N)(d\phi_t(e_i),d\phi_t(e_j))R^N(d\phi_t(e_i),d\phi_t(e_j))\tau(\phi_t)\big)d\phi_t(e_k)\\
&+2R^N\big(d\phi_t(e_k),R^N(\bar\nabla_{e_i}d\phi_t(\partial_t),d\phi_t(e_j))R^N(d\phi_t(e_i),d\phi_t(e_j))\tau(\phi_t)\big)d\phi_t(e_k)\\
&+R^N\big(d\phi_t(e_k),R^N(d\phi_t(e_i),d\phi_t(e_j))(\nabla_{d\phi_t(\partial_t)}R^N)(d\phi_t(e_i),d\phi_t(e_j))\tau(\phi_t)\big)d\phi_t(e_k)\\
&+2R^N\big(d\phi_t(e_k),R^N(d\phi_t(e_i),d\phi_t(e_j))R^N(\bar\nabla_{e_i}d\phi_t(\partial_t),d\phi_t(e_j))\tau(\phi_t)\big)d\phi_t(e_k)\\
&-R^N\big(d\phi_t(e_k),R^N(d\phi_t(e_i),d\phi_t(e_j))R^N(d\phi_t(e_i),d\phi_t(e_j))\bar\Delta d\phi_t(\partial_t)
\big)d\phi_t(e_k)\\
&+R^N\big(d\phi_t(e_k),R^N(d\phi_t(e_i),d\phi_t(e_j))R^N(d\phi_t(e_i),d\phi_t(e_j))
R^N(d\phi_t(\partial_t),d\phi_t(e_l))d\phi_t(e_l)\big)d\phi_t(e_k)\\
&+R^N\big(d\phi_t(e_k),R^N(d\phi_t(e_i),d\phi_t(e_j))R^N(d\phi_t(e_i),d\phi_t(e_j))\tau(\phi_t)\big)\bar\nabla_{e_k}d\phi_t(\partial_t).
\end{align*}

Evaluating at \(t=0\) then yields
\begin{align*}
\int_M&\langle d\phi_t(\partial_t),\frac{\bar\nabla}{\partial t}\big(\tr R^N(d\phi_t(\cdot),\Omega_0)d\phi_t(\cdot)\big)\rangle\dv\big|_{t=0}\\
=&
\int_M\langle V,(\nabla_{V}R^N)\big(d\phi(e_k),R^N(d\phi(e_i),d\phi(e_j))R^N(d\phi(e_i),d\phi(e_j))\tau(\phi)\big)d\phi(e_k)\rangle\dv\\
&+\int_M\langle V,R^N\big(\bar\nabla_{e_k}V,R^N(d\phi(e_i),d\phi(e_j))R^N(d\phi(e_i),d\phi(e_j))\tau(\phi)\big)d\phi(e_k)\rangle\dv\\
&+\int_M\langle V,R^N\big(d\phi(e_k),(\nabla_{V}R^N)(d\phi(e_i),d\phi(e_j))R^N(d\phi(e_i),d\phi(e_j))\tau(\phi)\big)d\phi(e_k)\rangle\dv\\
&+2\int_M\langle V,R^N\big(d\phi(e_k),R^N(\bar\nabla_{e_i}V,d\phi(e_j))R^N(d\phi(e_i),d\phi(e_j))\tau(\phi)\big)d\phi(e_k)\rangle\dv\\
&+\int_M\langle V,R^N\big(d\phi(e_k),R^N(d\phi(e_i),d\phi(e_j))(\nabla_{V}R^N)(d\phi(e_i),d\phi(e_j))\tau(\phi)\big)d\phi(e_k)\rangle\dv\\
&+2\int_M\langle V,R^N\big(d\phi(e_k),R^N(d\phi(e_i),d\phi(e_j))R^N(\bar\nabla_{e_i}V,d\phi(e_j))\tau(\phi)\big)d\phi(e_k)\rangle\dv\\
&-\int_M\langle V,R^N\big(d\phi(e_k),R^N(d\phi(e_i),d\phi(e_j))R^N(d\phi(e_i),d\phi(e_j))\bar\Delta 
V\big)d\phi(e_k)\rangle\dv\\
&+\int_M\langle V,R^N\big(d\phi(e_k),R^N(d\phi(e_i),d\phi(e_j))R^N(d\phi(e_i),d\phi(e_j))
R^N(V,d\phi(e_l))d\phi(e_l)\big)d\phi(e_k)\rangle\dv\\
&+\int_M\langle V,R^N\big(d\phi(e_k),R^N(d\phi(e_i),d\phi(e_j))R^N(d\phi(e_i),d\phi(e_j))\tau(\phi)\big)\bar\nabla_{e_k}V\rangle\dv.
\end{align*}
The claim now follows by adding up the various contributions.
\end{proof}
Making use of the previous Lemma we can now give our first result on the
stability properties of critical points of \(\hat E_4(\phi)\).

\begin{Satz}
\label{thm:harmonic-stable-escurv}
Let \(\phi\colon M\to N\) be a smooth harmonic map.
Then, it is a weakly stable critical point of the ES-4-curvature energy \(\hat E_4(\phi)\).
\end{Satz}
\begin{proof}
Inserting \(\tau(\phi)=0\) into \eqref{sv-es4-curvature} yields
\begin{align*}
&\frac{d^2}{dt^2}\frac{1}{4}\int_M|R^N(d\phi_t(e_i),d\phi_t(e_j))\tau(\phi_t)|^2\dv\big|_{t=0}\\
=&-\int_M\langle\bar\Delta V,R^N(d\phi(e_i),d\phi(e_j))R^N(d\phi(e_i),d\phi(e_j))\bar\Delta V\rangle\dv \\
&+\int_M\langle\bar\Delta V,R^N(d\phi(e_i),d\phi(e_j))R^N(d\phi(e_i),d\phi(e_j))R^N(V,d\phi(e_k))d\phi(e_k)\rangle\dv\\
&-\int_M\langle V,R^N\big(d\phi(e_k),R^N(d\phi(e_i),d\phi(e_j))R^N(d\phi(e_i),d\phi(e_j))\bar\Delta 
V\big)d\phi(e_k)\rangle\dv\\
&+\int_M\langle V,R^N\big(d\phi(e_k),R^N(d\phi(e_i),d\phi(e_j))R^N(d\phi(e_i),d\phi(e_j))
R^N(V,d\phi(e_l))d\phi(e_l)\big)d\phi(e_k)\rangle\dv.
\end{align*}
Exploiting the symmetries of the Riemann curvature tensor the above equation can be 
simplified to
\begin{align*}
\frac{d^2}{dt^2}\frac{1}{4}&\int_M|R^N(d\phi_t(e_i),d\phi_t(e_j))\tau(\phi_t)|^2\dv\big|_{t=0}\\
=&\int_M\big|R^N(d\phi(e_i),d\phi(e_j))\big(\bar\Delta V-R^N(V,d\phi(e_k))d\phi(e_k)\big)\big|^2\dv\\
\geq& 0,
\end{align*}
which already completes the proof.
\end{proof}

\subsection{The case of isometric immersions into space forms}
Now, we want to simplify the lengthy formula \eqref{sv-es4-curvature} by assuming
that the target manifold is a space form of constant curvature \(K\).
In addition, we will also assume that \(\phi\) is an isometric immersion
as we ultimately want to investigate the normal stability of the small ES-4-harmonic hypersphere
which is a map of this kind.

As \(\phi\) is an isometric immersion into a space form we get 
\begin{align*}
R^N(d\phi(X),d\phi(Y))\tau(\phi)=K\big(\langle d\phi(Y),\tau(\phi)\rangle d\phi(X)
-\langle d\phi(X),\tau(\phi)\rangle d\phi(Y)\big)
=0
\end{align*}
since \(\langle d\phi(Z),\tau(\phi)\rangle =0\) for all \(Z\in TM\).

Moreover, as we are considering a target with constant sectional curvature
all terms containing the derivative of the Riemann curvature tensor in \eqref{sv-es4-curvature}
will drop out.

\begin{Prop}
Let \(\phi\colon M\to N\) be a smooth map where \(N\) is a space form of constant curvature \(K\). 
Moreover, suppose that \(\langle d\phi(X),\tau(\phi)\rangle=0\) for all \(X\in TM\).
Then the second variation of \(\hat E^{ES}_4(\phi)\) simplifies to
\begin{align}
\label{sv-es4-spaceform-normal}
\frac{d^2}{dt^2}&\frac{1}{4}\int_M|R^N(d\phi_t(e_i),d\phi_t(e_j))\tau(\phi_t)|^2\dv\big|_{t=0} \\
\nonumber=&-\int_M\langle \frac{\bar\nabla}{\partial t}d\phi_t(\partial_t),\hat\tau_4(\phi_t)\rangle\big|_{t=0}\dv 
\\
\nonumber &
+K^2\int_M \bigg(
|\langle\bar\nabla V,\tau(\phi)\rangle|^2|d\phi|^2
-\langle\bar\nabla_{e_i}V,\tau(\phi)\rangle
\langle\bar\nabla_{e_j}V,\tau(\phi)\rangle\langle d\phi(e_i),d\phi(e_j)\rangle
\\
\nonumber &
+2\langle\bar\nabla_{e_i}V,\tau(\phi)\rangle\langle d\phi(e_i),d\phi(e_j)\rangle
\langle d\phi(e_j),\bar\Delta V\rangle
-2|d\phi|^2\langle d\phi(e_i),\bar\Delta V\rangle\langle\bar\nabla_{e_i}V,\tau(\phi)\rangle
\\
\nonumber &
+|d\phi|^2|\langle d\phi,\bar\Delta V\rangle|^2 
-\langle d\phi(e_i),d\phi(e_j)\rangle\langle d\phi(e_i),\bar\Delta V\rangle 
\langle d\phi(e_j),\bar\Delta V\rangle 
\bigg)\dv \\
&\nonumber+K^3\int_M\bigg(
-4|d\phi|^2\langle d\phi(e_j),V\rangle\langle d\phi(e_i),d\phi(e_j)\rangle
\langle\bar\nabla_{e_i}V,\tau(\phi)\rangle
+2|d\phi|^4\langle d\phi(e_i),V\rangle\langle\bar\nabla_{e_i}V,\tau(\phi)\rangle \\
\nonumber &+2\langle V,d\phi(e_k)\rangle\langle d\phi(e_j),d\phi(e_k)\rangle
\langle d\phi(e_i),d\phi(e_j)\rangle\langle\bar\nabla_{e_i}V,\tau(\phi)\rangle
-2|d\phi|^4\langle d\phi(e_j),V\rangle\langle\bar\Delta V,d\phi(e_j)\rangle 
\\ \nonumber &
-2\langle d\phi(e_j),\bar\Delta V\rangle\langle d\phi(e_i),d\phi(e_j)\rangle
\langle d\phi(e_i),d\phi(e_k)\rangle\langle V,d\phi(e_k)\rangle 
\bigg)\dv
\\ \nonumber &
+K^4\int_M\bigg(
|d\phi|^6|\langle V,d\phi\rangle|^2
-3|d\phi|^4\langle V,d\phi(e_i)\rangle\langle V,d\phi(e_j)\rangle\langle d\phi(e_i),d\phi(e_j)\rangle \\
\nonumber&+3|d\phi|^2\langle V,d\phi(e_i)\rangle\langle V,d\phi(e_k)\rangle
\langle d\phi(e_i),d\phi(e_j)\rangle \langle d\phi(e_j),d\phi(e_k)\rangle \\
\nonumber&-\langle V,d\phi(e_l)\rangle\langle V,d\phi(e_k)\rangle 
\langle d\phi(e_l),d\phi(e_j)\rangle \langle d\phi(e_j),d\phi(e_i)\rangle
\langle d\phi(e_k),d\phi(e_i)\rangle
\bigg)\dv.
\end{align}
\end{Prop}

\begin{proof}
By assumption, we have \(\langle d\phi(X),\tau(\phi)\rangle=0\) for all \(X\in TM\) such that
many of the terms in \eqref{sv-es4-curvature} vanish and we are left with
\begin{align}
&\frac{d^2}{dt^2}\frac{1}{4}
\int_M|R^N(d\phi_t(e_i),d\phi_t(e_j))\tau(\phi_t)|^2\dv\big|_{t=0} \\
\nonumber=&-\int_M\langle \frac{\bar\nabla}{\partial t}d\phi_t(\partial_t),\hat{\tau}_4(\phi_t)
\rangle\big|_{t=0}\dv \\
\nonumber&+\int_M\langle \bar\nabla_{e_i}V,R^N(R^N(\bar\nabla_{e_i}V,d\phi(e_j))\tau(\phi),\tau(\phi))d\phi(e_j)\rangle\dv
\\
\nonumber&+\int_M\langle \bar\nabla_{e_i}V,R^N(R^N(d\phi(e_i),\bar\nabla_{e_j}V)\tau(\phi),\tau(\phi))d\phi(e_j)\rangle\dv
\\
\nonumber&-\int_M\langle \bar\nabla_{e_i}V,R^N(R^N(d\phi(e_i),d\phi(e_j))\bar\Delta V,\tau(\phi))d\phi(e_j)
\rangle\dv
\\
\nonumber&+\int_M\langle \bar\nabla_{e_i}V,R^N(R^N(d\phi(e_i),d\phi(e_j))
R^N(V,d\phi(e_k))d\phi(e_k),\tau(\phi))d\phi(e_j)\rangle\dv
\\
\nonumber&+\int_M\langle\bar\Delta V,R^N(d\phi(e_i),d\phi(e_j))R^N(\bar\nabla_{e_i}V,d\phi(e_j))\tau(\phi)\rangle\dv \\
\nonumber&-\frac{1}{2}\int_M\langle\bar\Delta V,R^N(d\phi(e_i),d\phi(e_j))R^N(d\phi(e_i),d\phi(e_j))\bar\Delta V\rangle\dv \\
\nonumber&+\frac{1}{2}\int_M\langle\bar\Delta V,R^N(d\phi(e_i),d\phi(e_j))R^N(d\phi(e_i),d\phi(e_j))R^N(V,d\phi(e_k))d\phi(e_k)\rangle\dv\\
\nonumber&+\int_M\langle V,R^N\big(d\phi(e_k),R^N(d\phi(e_i),d\phi(e_j))R^N(\bar\nabla_{e_i}V,d\phi(e_j))\tau(\phi)\big)d\phi(e_k)\rangle\dv\\
\nonumber&-\frac{1}{2}\int_M\langle V,R^N\big(d\phi(e_k),R^N(d\phi(e_i),d\phi(e_j))R^N(d\phi(e_i),d\phi(e_j))\bar\Delta 
V\big)d\phi(e_k)\rangle\dv\\
\nonumber&+\frac{1}{2}\int_M\langle V,R^N\big(d\phi(e_k),R^N(d\phi(e_i),d\phi(e_j))R^N(d\phi(e_i),d\phi(e_j))
R^N(V,d\phi(e_l))d\phi(e_l)\big)d\phi(e_k)\rangle\dv.
\end{align}

In order to manipulate these remaining contributions we 
make use of the assumption that \(N\) is a space form of constant curvature such that the Riemann curvature tensor has the simple form \eqref{curvature-space-form}.
Now, also assuming that \(\langle d\phi(X),\tau(\phi)\rangle=0\) for all \(X\in TM\),  
we find 
\begin{align*}
\langle \bar\nabla_{e_i}V,R^N(R^N(\bar\nabla_{e_i}V,d\phi(e_j))\tau(\phi),\tau(\phi))d\phi(e_j)\rangle=&K^2|\langle\bar\nabla V,\tau(\phi)\rangle|^2|d\phi|^2,\\
\langle \bar\nabla_{e_i}V,R^N(R^N(d\phi(e_i),\bar\nabla_{e_j}V)\tau(\phi),\tau(\phi))d\phi(e_j)\rangle=&-K^2\langle\bar\nabla_{e_i}V,\tau(\phi)\rangle
\langle\bar\nabla_{e_j}V,\tau(\phi)\rangle\langle d\phi(e_i),d\phi(e_j)\rangle,\\
\langle \bar\nabla_{e_i}V,R^N(R^N(d\phi(e_i),d\phi(e_j))\bar\Delta V,\tau(\phi))d\phi(e_j)\rangle=&
K^2\big(
|d\phi|^2\langle d\phi(e_i),\bar\Delta V\rangle\langle\bar\nabla_{e_i} V,\tau(\phi)\rangle
 \\&
 -\langle\bar\nabla_{e_i}V,\tau(\phi)\rangle\langle d\phi(e_i),d\phi(e_j)\rangle
\langle d\phi(e_j),\bar\Delta V\rangle
\big),\\
|R^N(d\phi(e_i),d\phi(e_j))\bar\Delta V|^2=
&2K^2|d\phi|^2|\langle d\phi,\bar\Delta V\rangle |^2 \\
&-2K^2\langle d\phi(e_i),d\phi(e_j)\rangle\langle d\phi(e_i),\bar\Delta V\rangle
\langle d\phi(e_j),\bar\Delta V\rangle
\end{align*}
and also 
\begin{align*}
\langle&\bar\Delta V,R^N(d\phi(e_i),d\phi(e_j))R^N(\bar\nabla_{e_i}V,d\phi(e_j))\tau(\phi)\rangle \\
=&K^2\big(-|d\phi|^2\langle\bar\nabla_{e_i}V,\tau(\phi)\rangle
\langle\bar\Delta V,d\phi(e_i)\rangle
+\langle\bar\nabla_{e_i}V,\tau(\phi)\rangle
\langle d\phi(e_i),d\phi(e_j)\rangle \langle\bar\Delta V,d\phi(e_j)\rangle\big).
\end{align*}
Regarding the terms that contain the Riemann curvature tensor three times we find 
\begin{align*}
\langle& \bar\nabla_{e_i}V,R^N(R^N(d\phi(e_i),d\phi(e_j))R^N(V,d\phi(e_k))d\phi(e_k),\tau(\phi))d\phi(e_j)\rangle \\
=&K^3\big(-2|d\phi|^2\langle d\phi(e_j),V\rangle\langle d\phi(e_i),d\phi(e_j)\rangle
\langle\bar\nabla_{e_i}V,\tau(\phi)\rangle
+|d\phi|^4\langle d\phi(e_i),V\rangle\langle\bar\nabla_{e_i}V,\tau(\phi)\rangle \\
&+\langle V,d\phi(e_k)\rangle\langle d\phi(e_j),d\phi(e_k)\rangle
\langle d\phi(e_i),d\phi(e_j)\rangle\langle\bar\nabla_{e_i}V,\tau(\phi)\rangle\big).
\end{align*}
and 
\begin{align*}
\langle\bar\Delta V&,R^N(d\phi(e_i),d\phi(e_j))R^N(d\phi(e_i),d\phi(e_j))R^N(V,d\phi(e_k))d\phi(e_k)\rangle \\
=&K^3\big(-2|d\phi|^4\langle d\phi(e_j),V\rangle\langle\bar\Delta V,d\phi(e_j)\rangle \\
&-2\langle V,d\phi(e_k)\rangle\langle d\phi(e_j),d\phi(e_k)\rangle
\langle d\phi(e_i),d\phi(e_j)\rangle\langle\bar\Delta V,d\phi(e_i)\rangle \\
&+4|d\phi|^2\langle d\phi(e_i),V\rangle\langle\bar\Delta V,d\phi(e_j)\rangle
\langle d\phi(e_j),d\phi(e_i)\rangle
\big).
\end{align*}
Moreover, we obtain
\begin{align*}
\langle V&,R^N\big(d\phi(e_k),R^N(d\phi(e_i),d\phi(e_j))R^N(\bar\nabla_{e_i}V,d\phi(e_j))\tau(\phi)\big)d\phi(e_k)\rangle \\
=&K^3\big(-|d\phi|^2\langle\bar\nabla_{e_i}V,\tau(\phi)\rangle
\langle d\phi(e_i),d\phi(e_k)\rangle\langle V,d\phi(e_k)\rangle
+|d\phi|^4\langle\bar\nabla_{e_i}V,\tau(\phi)\rangle\langle V,d\phi(e_i)\rangle \\
&+\langle\bar\nabla_{e_i}V,\tau(\phi)\rangle
\langle d\phi(e_i),d\phi(e_j)\rangle\langle d\phi(e_j),d\phi(e_k)\rangle
\langle V,d\phi(e_k)\rangle \\
&-|d\phi|^2\langle\bar\nabla_{e_i}V,\tau(\phi)\rangle\langle d\phi(e_i),d\phi(e_j)\rangle
\langle V,d\phi(e_j)\rangle\big)
\end{align*}
and
\begin{align*}
\langle V&,R^N\big(d\phi(e_k),R^N(d\phi(e_i),d\phi(e_j))R^N(d\phi(e_i),d\phi(e_j))\bar\Delta 
V\big)d\phi(e_k)\rangle \\
=&K^3\big(2\langle d\phi(e_j),\bar\Delta V\rangle\langle d\phi(e_i),d\phi(e_j)\rangle
\langle d\phi(e_i),d\phi(e_k)\rangle\langle V,d\phi(e_k)\rangle \\
&-4|d\phi|^2\langle d\phi(e_j),\bar\Delta V\rangle\langle d\phi(e_i),d\phi(e_j)\rangle
\langle V,d\phi(e_i)\rangle
+2|d\phi|^4\langle d\phi(e_j),\bar\Delta V\rangle\langle V,d\phi(e_j)\rangle\big).
\end{align*}
Finally, for the term being quartic in the curvature tensor we find
\begin{align*}
\langle &V,R^N\big(d\phi(e_k),R^N(d\phi(e_i),d\phi(e_j))R^N(d\phi(e_i),d\phi(e_j))
R^N(V,d\phi(e_l))d\phi(e_l)\big)d\phi(e_k)\rangle \\
=&K^4\big(2|d\phi|^6|\langle V,d\phi\rangle|^2
-6|d\phi|^4\langle V,d\phi(e_i)\rangle\langle V,d\phi(e_j)\rangle\langle d\phi(e_i),d\phi(e_j)\rangle \\
&+6|d\phi|^2\langle V,d\phi(e_i)\rangle\langle V,d\phi(e_k)\rangle
\langle d\phi(e_i),d\phi(e_j)\rangle \langle d\phi(e_j),d\phi(e_k)\rangle \\
&-2\langle V,d\phi(e_l)\rangle\langle V,d\phi(e_k)\rangle 
\langle d\phi(e_l),d\phi(e_j)\rangle \langle d\phi(e_j),d\phi(e_i)\rangle
\langle d\phi(e_k),d\phi(e_i)\rangle\big).
\end{align*}
The claim follows by combining all the contributions.
\end{proof}
A direct consequence of the previous calculation is the following
\begin{Cor}
Let \(\phi\colon M\to N\) be a smooth map where \(N\) is a space form of constant curvature \(K\).
Moreover, suppose that \(\langle d\phi(X),\tau(\phi)\rangle=0\) for all \(X\in TM\)
and that \(\langle V,d\phi(Y)\rangle=0\) for all \(Y\in TM\). Then, the second variation
of the energy functional \(\hat E_4(\phi)\) simplifies to
\begin{align}
\label{sv-es4-spaceform-normal-simplified}
\frac{d^2}{dt^2}&\frac{1}{4}\int_M|R^N(d\phi_t(e_i),d\phi_t(e_j))\tau(\phi_t)|^2\dv\big|_{t=0} \\
\nonumber=&-\int_M\langle \frac{\bar\nabla}{\partial t}d\phi_t(\partial_t),\hat\tau_4(\phi_t)\rangle\big|_{t=0}\dv 
\\
\nonumber &
+K^2\int_M \bigg(
|\langle\bar\nabla V,\tau(\phi)\rangle|^2|d\phi|^2
-\langle\bar\nabla_{e_i}V,\tau(\phi)\rangle
\langle\bar\nabla_{e_j}V,\tau(\phi)\rangle\langle d\phi(e_i),d\phi(e_j)\rangle
\\
\nonumber &
+2\langle\bar\nabla_{e_i}V,\tau(\phi)\rangle\langle d\phi(e_i),d\phi(e_j)\rangle
\langle d\phi(e_j),\bar\Delta V\rangle
-2|d\phi|^2\langle d\phi(e_i),\bar\Delta V\rangle\langle\bar\nabla_{e_i}V,\tau(\phi)\rangle
\\
\nonumber &
+|d\phi|^2|\langle d\phi,\bar\Delta V\rangle|^2 
-\langle d\phi(e_i),d\phi(e_j)\rangle\langle d\phi(e_i),\bar\Delta V\rangle 
\langle d\phi(e_j),\bar\Delta V\rangle 
\bigg)\dv.
\end{align}
\end{Cor}
\begin{proof}
This follows directly from \eqref{sv-es4-spaceform-normal} taking into account \(\langle V,d\phi(Y)\rangle =0\) for all \(Y\in TM\).
\end{proof}

\subsection{Application to ES-4-harmonic hyperspheres}
In this subsection we will apply the general formula for the second variation of \(\hat E_4(\phi)\) to the case of the small hypersphere given by \(\phi\colon\s^m(\frac{1}{2})\to\s^{m+1}\).
Employing \eqref{eq:laplace-hypersurface} we now find:

\begin{Prop}
Let \(\phi\colon\s^m(\frac{1}{2})\to\s^{m+1}\) be the small hypersphere considered
as a critical point of \(\hat E_4(\phi)\). Then the quadratic form describing its 
normal index is given by
\begin{align}
\label{quadratic-form-es4-hyper-a}
\hat Q^{ES}_4(f\nu,f\nu)=(m-1)\int_{\s^m(\frac{1}{2})}&\big(
m^2H^2|\nabla f|^2-4mH\langle A(\nabla f),\nabla f\rangle+4|A(\nabla f)|^2
\big)\dv,
\end{align}
where \(\nu\) is the unit normal of the hypersphere \(\s^m(\frac{1}{2})\) and \(f\in C_c^\infty(\s^m(\frac{1}{2}))\).
\end{Prop}

\begin{proof}
We use the general formula \eqref{sv-es4-spaceform-normal-simplified}, set \(K=1\),
and simplify the terms as follows
\begin{align*}
|\langle\bar\nabla (f\nu),\tau(\phi)\rangle|^2|d\phi|^2=&m^3H^2|\nabla f|^2,&
\\
\langle\bar\nabla_{e_i}(f\nu),\tau(\phi)\rangle
\langle\bar\nabla_{e_j}(f\nu),\tau(\phi)\rangle\langle d\phi(e_i),d\phi(e_j)\rangle
=&m^2H^2|\nabla f|^2
,\\
\langle\bar\nabla_{e_i}(f\nu),\tau(\phi)\rangle\langle d\phi(e_i),d\phi(e_j)\rangle
\langle d\phi(e_j),\bar\Delta (f\nu)\rangle =&2mH\langle\nabla f,A(\nabla f)\rangle
,\\
|d\phi|^2\langle d\phi(e_i),\bar\Delta(f\nu)\rangle\langle\bar\nabla_{e_i}(f\nu),\tau(\phi)\rangle
=&2m^2H\langle A(\nabla f),\nabla f\rangle
,\\
|d\phi|^2|\langle d\phi,\bar\Delta (f\nu)\rangle |^2=&
4m|A(\nabla f)|^2
,\\
\langle d\phi(e_i),d\phi(e_j)\rangle\langle d\phi(e_i),\bar\Delta (f\nu)\rangle 
\langle d\phi(e_j),\bar\Delta (f\nu)\rangle 
=&4|A(\nabla f)|^2,
\end{align*}
which already completes the proof.
\end{proof}

Now, we can easily obtain the following result.

\begin{Satz}
\label{thm:hypersphere-crit-curv}
Let \(\phi\colon\s^m(\frac{1}{2})\to\s^{m+1}\) be the small hypersphere considered
as a critical point of \(\hat E^{ES}_4(\phi)\). Then the quadratic form describing its 
normal index is given by
\begin{align}
\label{quadratic-form-es4-hyper-b}
\hat Q^{ES}_4(f\nu,f\nu)=3(m-1)(m-2)^2\int_{\s^m(\frac{1}{2})}
|\nabla f|^2
\dv.
\end{align}
In particular, the small hypersphere \(\phi\colon\s^m(\frac{1}{2})\to\s^{m+1}\)
is a weakly stable critical point of \(\hat E_4(\phi)\) with respect to normal variations.
\end{Satz}

\begin{proof}
Recall that for the small hypersphere we have
\begin{align*}
A=-\sqrt{3}\operatorname{Id},\qquad H=-\sqrt{3},\qquad |A|^2=3m,\qquad H^2=3.
\end{align*}
Inserting this geometric data into \eqref{quadratic-form-es4-hyper-a} directly
yields the result.
\end{proof}

\begin{Bem}
Note that \eqref{quadratic-form-es4-hyper-b} 
vanishes for \(m=1\) which is to be expected as the curvature term
in \(E^{ES}_4(\phi)\) is zero in the case of a one-dimensional domain.
However, for \(\dim M\geq 3\) the quadratic form \eqref{quadratic-form-es4-hyper-b} will always be
non-negative and thus not contribute to the normal index.
\end{Bem}

\subsection{On the normal stability of the small ES-4-harmonic hypersphere}
Finally, we will now combine the results obtained on the stability of
critical points of both \(E_4(\phi)\) and \(\hat E_4(\phi)\) such that we
can make statements on the normal stability of \(ES-4\)-harmonic maps.

Recall that
\begin{align*}
E^{ES}_4(\phi)=E_4(\phi)+\hat E_4(\phi)
\end{align*}
with the critical points
\begin{align*}
\tau^{ES}_4(\phi):=\tau_4(\phi)+\hat\tau_4(\phi)=0,
\end{align*}
see \eqref{es-4-tension} for the precise form of the \(ES-4\)-tension field.

Our first, more general, result is the following:

\begin{Satz}
\label{thm:harmonic-stable-es4}
A smooth harmonic map \(\phi\colon M\to N\) is a weakly stable critical point of \(E^{ES}_4(\phi)\).
\end{Satz}
\begin{proof}
Using the assumption \(\tau(\phi)=0\) in both \eqref{sv-4energy-simplified} and \eqref{sv-es4-curvature} we find
\begin{align*}
\frac{d^2}{dt^2}E^{ES}_4(\phi_t)\big|_{t=0}
=&\frac{d^2}{dt^2}E_4(\phi_t)\big|_{t=0}+\frac{d^2}{dt^2}\hat E_4(\phi_t)\big|_{t=0}\\
=&\int_M\big|\bar\Delta^2V-\bar\Delta\big(R^N(V,d\phi(e_k))d\phi(e_k)\big)\big|^2\dv  \\
&+
\int_M\big|R^N(d\phi(e_i),d\phi(e_j))(\bar\Delta V-R^N(V,d\phi(e_k))d\phi(e_k)\big)\big|^2\dv\\
\geq& 0
\end{align*}
completing the proof.
\end{proof}

As the \(ES-4\) energy is the sum of the \(4\)-energy and the \(ES-4\) curvature term
the quadratic form describing the normal stability of an \(ES-4\)-harmonic map also
is the sum of the two corresponding quadratic forms, i.e.

\begin{align}
Q^{ES}_{4}(f\nu,f\nu):=Q_4(f\nu,f\nu)+\hat{Q}_4(f\nu,f\nu).
\end{align}

We are now ready to prove the main result of this article.

\begin{Satz}
\label{thm:stability-es4}
Let \(\phi\colon\s^m(\frac{1}{2})\to\s^{m+1}\)
be the small proper \(ES-4\)-harmonic hypersphere. Then, the quadratic form 
describing its normal stability is given by
\begin{align}
\label{qf4es}
Q^{ES}_{4}(f\nu,f\nu)=\int_{\s^m(\frac{1}{2})}
\big(&|\Delta^2f|^2+(10m+72)|\nabla\Delta f|^2\\
\nonumber&+(25m^2+286m+480)|\Delta f|^2 \\
\nonumber&+(-33m^3+312m^2-900m+576)|\nabla f|^2 \\
\nonumber&-216m^4f^2
\big)\dv.
\end{align}
In particular, the normal index of the small proper ES-4-harmonic hypersphere
is equal to one, i.e.
\begin{align*}
\operatorname{Ind}^{\rm{nor}}_{ES-4}(\s^m(\frac{1}{2})\to\s^{m+1})=1.
\end{align*}
\end{Satz}

\begin{proof}
The first claim of the proof follows directly from adding up the corresponding quadratic forms \eqref{qf4-c}
and \eqref{quadratic-form-es4-hyper-b}.

Regarding the second claim we assume that \(\lambda\) is an eigenvalue of the Laplacian, that is \(\Delta f=\lambda f\) such that the quadratic form \eqref{qf4es} becomes
\begin{align*}
Q^{ES}_{4}(f\nu,f\nu)=\int_{\s^m(\frac{1}{2})}
\big(&\lambda^4+(10m+72)\lambda^3\\
\nonumber&+(25m^2+286m+480)\lambda^2 \\
\nonumber&+(-33m^3+312m^2-900m+576)\lambda \\
\nonumber&-216m^4f^2
\big)\dv.
\end{align*}
For \(\lambda=0\) corresponding to the constant eigenfunction of the Laplacian 
the above quadratic form is negative producing the normal index of the hypersphere.

Now, recall that the first non-zero eigenvalue of the Laplace operator on \(\s^m(\frac{1}{2})\) 
is \(\lambda_1=4m\) in which case the above quadratic form becomes
\begin{align*}
Q^{ES}_{4}(f\nu,f\nu)=\int_{\s^m(\frac{1}{2})}
\big(
948m^4+10432m^3+4080m^2+2304m
\big)f^2\dv.
\end{align*}
This is clearly positive for all values of \(m\) and all higher eigenvalues of the Laplace operator
concluding the proof.
\end{proof}

\bibliographystyle{plain}
\bibliography{mybib}
\end{document}